\documentclass[11pt, reqno]{article}
\usepackage{}
\usepackage[total={6in, 9in}]{geometry}
\usepackage{amsthm}
\usepackage{amssymb}
\usepackage{amsmath,amsfonts}
\usepackage{mathrsfs}
\usepackage{cases}
\usepackage{latexsym,bm}
\usepackage{indentfirst}
\usepackage{color}
\usepackage{ifpdf}
\usepackage{graphicx}
\usepackage{psfrag}
\usepackage{dsfont}
\usepackage[ruled,vlined]{algorithm2e}
\usepackage{enumerate}
\usepackage[
pdfauthor={},
pdftitle={Edge coloring graphs with large minimum degree II},
pdfstartview=XYZ,
bookmarks=true,
colorlinks=true,
linkcolor=blue,
urlcolor=blue,
citecolor=blue,
bookmarks=true,
linktocpage=true,
hyperindex=true
]{hyperref}

\usepackage{pgf,tikz,pgfplots}
\usepackage{tikz-3dplot}
\usepackage{pifont}
\usepackage{refcount}
\newtheorem{THM}{\textbf{Theorem}}[section]

\newtheorem{LEM}[THM]{\textbf{Lemma}}

\newtheorem{CON}[THM]{\textbf{Conjecture}}

\newcommand{\pf}{\noindent\textbf{Proof}.\quad}

\newtheorem*{THM1}{\textbf{Theorem 1.3}}
\newtheorem*{THM2}{\textbf{Theorem 1.4}}

\newcommand{\ve}{\varepsilon}

\newcommand{\CC}{\mathcal{C}}

\DeclareMathOperator{\df}{def}
\newcommand{\pbar}{\overline{\varphi}}

\newcommand{\arxiv}[1]{\href{http://arxiv.org/abs/#1}{\texttt{arXiv:#1}}}

\begin{document}
\title{The overfull conjecture on graphs of odd order and large minimum degree}
\author	{Songling Shan \\ 
	\medskip  Illinois State  University, Normal, IL 61790\\
	\medskip 
	{\tt sshan12@ilstu.edu}
}

\date{\today}
\maketitle

\emph{\textbf{Abstract}.}
Let $G$ be a simple graph with maximum degree $\Delta(G)$. A subgraph $H$ of $G$ is overfull if $|E(H)|>\Delta(G)\lfloor \frac{1}{2}|V(H)| \rfloor$. Chetwynd and Hilton in 1986 conjectured that a graph $G$  with $\Delta(G)>\frac{1}{3}|V(G)|$
has chromatic index $\Delta(G)$ if and only if $G$ contains no overfull subgraph. 
 Let $0<\ve <1$ and  $G$ be a large graph on $n$ vertices with minimum degree at least $\frac{1}{2}(1+\ve)n$. It was shown 
that the conjecture holds for $G$ if $n$ is even. In this paper, 
the same result is proved  if $n$ is odd.  As far as we know, this is the first result on
the conjecture for graphs of odd order and with a minimum degree constraint.


\vspace{2mm}

\section{Introduction}

In this paper,  the terminology  graph is used to represent a simple graph. 
A multigraph may contain multiple edges but contains no loop. 
Let $G$ be a multigraph.
We use $V(G)$ and  $E(G)$ to denote the vertex set and the edge set of $G$,
respectively, and  let  $e(G)=|E(G)|$. 
For $v\in V(G)$, $N_G(v)$ is the set of neighbors of $v$ 
in $G$, and 
$d_G(v)$, the degree of $v$
in $G$, is the number of edges of $G$ that are incident with $v$.
We let $d_G^s(v)=|N_G(v)|$, and call it the \emph{simple degree}
of $v$ in $G$.   
When $G$ is simple, we have $d_G(v)=d_G^s(v)$.  
Let $V_1,
V_2\subseteq V(G)$ be two disjoint vertex sets. Then $E_G(V_1,V_2)$ is the set
of edges in $G$  with one end in $V_1$ and the other end in $V_2$, and  $e_G(V_1,V_2):=|E_G(V_1,V_2)|$.  We write $E_G(v,V_2)$ and $e_G(v,V_2)$
if $V_1=\{v\}$ is a singleton.  
We also use $G[V_1,V_2]$
to denote the bipartite subgraph of $G$ with vertex set $V_1\cup V_2$
and edge set $E_G(V_1, V_2)$. 
For 
$S\subseteq V(G)$ and $v\in V(G)$, $N_G(v, S):=N_G(v)\cap S$ and $d_S(v):=e_G(v,S\setminus\{v\})$,  
 the subgraph of $G$ induced by  $S$ is  $G[S]$, and  $G-S:=G[V(G)\setminus S]$. 
 If $F\subseteq E(G)$, then $G-F$ is obtained from $G$ by deleting all
 the edges of $F$. 
  Let $\mu(G)=\max\{e_G(u,v)\,:\, u,v\in V(G)\}$ be 
the multiplicity of $G$, and  for $x\in V(G)$, let $\mu(x)=\max\{e_G(x,v): v\in V(G)\}$
be the multiplicity of $x$. 

For two integers $p,q$, let $[p,q]=\{ i\in \mathbb{Z} \,:\, p \le i \le q\}$.
For a nonnegative integer $k$,  
 an  \emph{edge $k$-coloring} of a multigraph $G$ is a mapping $\varphi$ from $E(G)$ to the set of integers
$[1,k]$, called  \emph{colors}, such that  no two adjacent edges receive the same color with respect to $\varphi$.  
The \emph{chromatic index} of $G$, denoted $\chi'(G)$, is  the smallest integer $k$ so that $G$ has an edge $k$-coloring.  
We denote by $\CC^k(G)$ the set of all edge $k$-colorings of $G$.

In 1960's, Gupta~\cite{Gupta-67}  and, independently, Vizing~\cite{Vizing-2-classes}  showed
 that for all graphs $G$,  $\Delta(G) \le \chi'(G) \le \Delta(G)+1$. 
This 
leads to a natural classification of   graphs. Following Fiorini and Wilson~\cite{fw},   a graph $G$ is of \emph{class 1} if $\chi'(G) = \Delta(G)$ and of \emph{class 2} if $\chi'(G) = \Delta(G)+1$.  Holyer~\cite{Holyer} showed that it is NP-complete to determine whether an arbitrary graph is of class 1. 
However, there are sufficient conditions for a graph to be class 2. 
For example,  any graph $G$ with $|E(G)|>\Delta(G) \lfloor |V(G)|/2\rfloor$ is class 2, 
since each color class of $G$ is a matching and each matching in $G$
has size at most $\lfloor |V(G)|/2\rfloor$, 
we have 
$\chi'(G) \ge  |E(G)|/ \lfloor |V(G)|/2\rfloor>\Delta(G)$ but $\chi'(G) \le \Delta(G)+1$ by the result of Gupta and Vizing. 
Such graphs $G$ are  called \emph{overfull} (we can define a multigraph to be overfull in the same way).
An overfull subgraph $H$ of $G$ with  $\Delta(H)=\Delta(G)$
is called a \emph{$\Delta(G)$-overfull subgraph} of $G$. 
 A number of long-standing conjectures listed in {\it Twenty Pretty Edge Coloring Conjectures} in~\cite{StiebSTF-Book} lie in deciding when a 
 graph is overfull.   Chetwynd and  Hilton~\cite{MR848854,MR975994},  in 1986, proposed the following 
conjecture. 
\begin{CON}[Overfull Conjecture]\label{overfull-con}
	Let $G$ be a  graph  with $\Delta(G)>\frac{1}{3}|V(G)|$. Then $\chi'(G)=\Delta(G)$  if and only if $G$ contains no $\Delta(G)$-overfull subgraph.  
\end{CON}

The $3$-critical graph $P^*$, obtained from the Petersen graph by deleting one vertex, has $\chi'(P^*)=4$, 
satisfies $\Delta(P^*)=\frac{1}{3}|V(P^*)|$ but contains no  3-overfull subgraph. 
Thus the degree condition  $\Delta(G)>\frac{1}{3}|V(G)|$ in the conjecture above is best possible.  Applying Edmonds' matching polytope theorem, Seymour~\cite{seymour79}  showed  that whether a graph  $G$ contains an overfull subgraph of maximum degree $\Delta(G)$ can be determined in polynomial time. Thus if the Overfull Conjecture is true, then the NP-complete problem of 
determining the chromatic index becomes  polynomial-time solvable 
for graphs $G$ with $\Delta(G)>\frac{1}{3}|V(G)|$.
There have been some fairly strong results supporting the Overfull Conjecture in the case when $G$ is regular.  It is easy to verify that when $G$ is regular with even order, $G$ has no $\Delta(G)$-overfull subgraphs if its vertex degrees are at least $\frac{1}{2}|V(G)|$.  Thus the well-known 1-Factorization Conjecture stated below is a special case of the Overfull Conjecture.

\begin{CON}[1-Factorization Conjecture]
	Let  $G$ be a regular graph of order  $2n$  with degree at least $n$ if  $n$ is odd, or at least $n-1$ if $n$ is even.  Then $G$ is 1-factorable; equivalently, $ \chi'(G) = \Delta(G)$.
\end{CON}

Hilton and Chetwynd~\cite{MR1001390} verified the 1-Factorization Conjecture if the vertex degree is at least  $0.823|V(G)|$.  Perkovi\'c and Reed~\cite{MR1439301} showed in 1997 that the 1-Factorization Conjecture is true for large regular graphs with vertex degree at least $|V(G)|/(2-\ve)$ for any $0<\ve<1$.   In 2016,  Csaba, K\"uhn, Lo, Osthus and Treglown~\cite{MR3545109} verified the conjecture for sufficiently large $|V(G)|$. 
Much less is known about the truth of the Overfull Conjecture if we no longer require that $G$ is regular.  It was confirmed for graphs with $\Delta(G) \ge |V(G)| -3$ by Chetwynd and Hilton in 1989~\cite{MR975994}. Plantholt~\cite{MR2082738} in 2004 verified the conjecture for graphs of even order and minimum degree at least  $0.8819|V(G)|$.  More recently, Plantholt~\cite{Plantholt2} showed the conjecture  for sufficiently large even order  graphs with minimum degree at least $\frac{2}{3}|V(G)|$.  Plantholt~\cite{2104.06253} and the author extended these results and gave an asymptotic result for general graphs that is similar to the Perkovi\'c-Reed result for regular graphs, by showing the following result: 
	For all $0<\ve <1$, there exists $n_0$
	such that the following statement holds:
	if $G$ is a graph on $2n\ge n_0$ vertices with $\delta(G) \ge (1+\ve)n$, then $\chi'(G)=\Delta(G)$ if and only if $G$ contains no $\Delta(G)$-overfull subgraph. Furthermore, there is a polynomial time algorithm that finds an optimal coloring. 
	All the results mentioned above on graphs with minimum degree constraints are for graphs of 
	even order. 
 In this paper, we obtain the first  result (as far as we know) on the Overfull Conjecture with large minimum degree constraints for graphs of odd order.

\begin{THM}\label{thm:1}
For all $0<\ve <1$, there exists $n_0$
such that the following statement holds:
if $G$ is a graph on $2n-1\ge n_0$ vertices with $\delta(G) \ge (1+\ve)n$,   then $\chi'(G)=\Delta(G)$ if and only if $G$  is not overfull. Furthermore, there is a polynomial time algorithm that finds an optimal coloring. 
\end{THM}

The remainder of this paper is organized as follows.
In the next section, we introduce some notation and   preliminary results.
In   Section 3, we prove 
Theorem~\ref{thm:1} by applying Theorem~\ref{thm:D-coloring} from Section 3. 
Theorem~\ref{thm:D-coloring}
is then proved in the last section.

\section{Notation and preliminaries}

\subsection{Results on degree sequences, paths, matchings, and overfull graphs}

We need the following classic result of Hakimi~\cite{MR148049} on multigraphic degree sequence. 
\begin{THM}\label{thm:degree-seq}
	Let $0 \le d_n \le \ldots \le d_1$ be integers. Then there exists a multigraph $G$
	on vertices $x_1,\ldots, x_n$ such that $d_G(x_i)=d_i$
	for all $i$ if and only if $\sum_{i=1}^nd_i$ is even and $\sum_{i>1}d_i \ge d_1$. 
\end{THM}

Though it is not explicitly stated in~\cite{MR148049}, the inductive proof yields a polynomial
time algorithm which finds an appropriate multigraph if it exists.

\begin{THM}[\cite{MR47308}]\label{thm:Dirac}
	If  $G$ is a graph on $n\ge 3$ vertices with $\delta(G) \ge \frac{n}{2}$, then $G$ is hamiltonian.
\end{THM}

Following the proof of Dirac~\cite{MR47308}, a hamiltonian cycle can be constructed in polynomial time in $n$
if $\delta(G) \ge \frac{n}{2}$. In fact, there is a polynomial time algorithm that constructs the closure of a graph $G$ and finds a hamiltonian cycle of $G$
if its closure is a complete graph (see~\cite[Exercise 4.2.15, page 62]{MR0411988}). 

A path $P$ connecting two vertices $u$ and $v$ is called 
a {\it $(u,v)$-path}, and we write $uPv$ or $vPu$ in order to specify the two endvertices of 
$P$. Let $uPv$ and $xQy$ be two disjoint paths. If $vx$ is an edge, 
we write $uPvxQy$ as
the concatenation of $P$ and $Q$ through the edge $vx$. If $P$ is a path and $x,y\in V(P)$,
then $xPy$ is the subpath of $P$ with endvertices $x$ and $y$.  

We will use the following notation: $0<a \ll b \le 1$. 
Precisely, if we say a claim is true provided that $0<a \ll b \le 1$, 
then this means that there exists a non-decreasing function $f:(0,1]\rightarrow (0,1]$ such that the statement holds for all $0<a,b\le 1$ satisfying $a \le f(b)$. 

\begin{LEM}[\cite{2104.06253}]\label{lem:path-decomposition0}
	Let $0<\frac{1}{n_0} \ll  \ve <1$, and $G$ be 
	graph on $n\ge n_0$ vertices such that $\delta(G) \ge \frac{1}{2}(1+\ve) n$. 
	Moreover, let $M=\{a_1b_1,\ldots,a_tb_t\}$ be a matching in the complete graph on $V(G)$ of size at most $\frac{\ve n}{8}$. 
	Then there exist vertex-disjoint paths $P_1,\ldots, P_t$ in $G$ such that $\bigcup V(P_i)=V(G)$
	and $P_i$ joins $a_i$ to $b_i$,  and these paths can be found in polynomial time.
\end{LEM}

A multigraph $G$ is a \emph{star-multigraph} if $G$
has a vertex $x$ that is incident with all multiple edges of $G$. 
In other words, $G-x$ is a simple graph.  The vertex $x$ is called the \emph{multi-center}
of $G$.

\begin{LEM}\label{lem:path-decomposition}
	Let $0<\frac{1}{n_0} \ll  \eta\ll \ve  <1$, and $G$ be start-multigraph 
	graph on $n\ge n_0$ vertices such that $\delta(G) \ge \frac{1}{2}(1+\ve) n$ and $ \mu(x) <  \eta n $,  where $x$ is the multi-center of $G$.  
	Moreover, let $M=\{a_1b_1,\ldots,a_tb_t\}$ be a matching in the complete graph on $V(G)$ of size at most $\frac{\ve n}{13}$. If $|N_G(x)\setminus \{a_1,b_1,\ldots, a_t,b_t\}|\ge 2$,  
	then there exist vertex-disjoint paths $P_1,\ldots, P_t$ in $G$ such that $\bigcup V(P_i)=V(G)$
	and $P_i$ joins $a_i$ to $b_i$,  and these paths can be found in polynomial time. 
\end{LEM}
\pf By relabeling the matching edges if necessary,  assume that 
if $x\in \{a_1,b_1, \ldots, a_t,b_t\}$, then $x=a_t$.  

If $x\in \{a_1,b_1,\ldots, a_t,b_t\}$ and so $x=a_t$,    
then let $a_t'\in N_G(x)\setminus \{a_1,b_1,\ldots, a_t,b_t\}$,
and $M'=(M\setminus \{a_tb_t\} )\cup \{a_t'b_t\}$. 
If $x\not\in \{a_1,b_1, \ldots, a_t,b_t\}$, then let $x_1, x_2\in N_G(x)\setminus \{a_1,b_1, \ldots, a_t,b_t\}$ be distinct, 
and let $M'=(M\setminus \{a_tb_t\} )\cup \{a_tx_1, x_2b_t\}$. 

Note that $\delta(G-x) \ge \frac{1}{2}(1+\ve) n-\eta n \ge\frac{1}{2}(1+\frac{2\ve}{3})n$
and $|M'| \le  \frac{\ve n}{13}+1<   \frac{\ve n}{12}$. 
Applying 
Lemma~\ref{lem:path-decomposition0} to $G-x$ with matching $M'$, 
we find vertex-disjoint paths $P'_1,\ldots,  P'_t$ in $G-x$ such that $\bigcup V(P'_i)=V(G-x)$. Furthermore, 
if $x=a_t$, 
$P'_i$ joins $a_i$ to $b_i$ for $i\in [1,t-1]$, and $P_t'$ joins $a_t'$
to $b_t$; if $x \not\in \{a_1,b_1, \ldots, a_t,b_t\}$, 
$P'_i$ joins $a_i$ to $b_i$ for $i\in [1,t-1]$,  $P'_{t}$ joins $a_t$ to $x_1$, and $P_{t+1}'$ joins $x_2$
to $b_t$.   
Letting $P_i=P_i'$ for $i\in [1,t-1]$,  and  $P_t=a_ta_t'P_t'b_t$ if $x=a_t$ 
and $ P_t=a_tP_t'x_1xx_2P'_{t+1}b_t$ if $x \not\in \{a_1,b_1, \ldots, a_t,b_t\}$ 
gives the desired paths for $G$ and $M$. By  Lemma~\ref{lem:path-decomposition0} and the simple adjustment on $M$ above,
it is clear that these paths can be found also in polynomial time. 
\qed

\begin{LEM}\label{thm:matching3}
	Let $G$ be a   graph of order $2n$ and minimum degree at least 1. Assume that except at most one vertex all other vertices have degree at least $n+1$ in $G$.    Then $G$ has a perfect matching, and the matching can be found in polynomial time.  
\end{LEM}

\pf Let $u \in V(G)$ such that $d_G(u)=\delta(G)$. Since $\delta(G)\ge 1$, we let $v$ be a neighbor of $u$, and let $G_1=G-\{u,v\}$.
Since $\delta(G-u)\ge n$, we have   
$\delta(G_1) \ge n-1= \frac{1}{2}(2n-2)$. Applying Theorem~\ref{thm:Dirac}, $G_1$
has a perfect matching $M$. Thus $M\cup \{uv\}$ is a perfect matching of $G$.
By the construction of $G_1$ and the comments immediately below  Theorem~\ref{thm:Dirac}, 
such a perfect matching of $G$ can be found in polynomial time. 
\qed

\begin{LEM}\label{lem:matching-in-bipartite}
	Let $G[X,Y]$ be bipartite star-multigraph with $|X|=|Y|=n$ and multi-center $x$. Suppose that $\delta(G)\ge 1$ and $\delta(G-x) \ge t$ 
	for some $t\in [1,n-1]$, 
	and except at most $t$ vertices all other vertices of $G$ have simple degree 
	at least $\frac{n}{2}+1$ in $G$. Then 
	$G$ has a perfect matching. Furthermore, such a matching can be found in polynomial time. 
\end{LEM}
\pf Suppose, without loss of generality, that $x\in X$. 
As $\delta(G) \ge 1$, we let $y\in Y$ be a neighbor of $x$. Let $X'=X\setminus\{x\}$
and $Y'=Y\setminus\{y\}$. 
It suffices to show that $G':=G[X', Y']$ satisfies Hall's Condition. If not, 
we let $S\subseteq X'$ with smallest cardinality such that 
$|S|>|N_{G'}(S)|$.  Since $\delta(G')\ge 1$, it follows  that $|S|\ge 1$. 
Then by the choice of $S$,  we have $|S|=|N_{G'}(S)|+1$ and $|N_{G'}(S)|<|Y'|$.  
As  $|S|>|N_{G'}(S)|$, it follows that 
$|S| \ge \delta(G')+1 \ge t+1$. As $G'$
has at most $t$ vertices of  degree less than $\frac{n}{2}$, $S$ contains a vertex 
of degree at least $\frac{n}{2}$ in $G'$ and so $|N_{G'}(S)|\ge \frac{n}{2}$. 
Thus $|S|>\frac{n}{2}$ and consequently  
 $|X'\setminus S|  <\frac{n}{2}$. Since $|N_{G'}(S)|<|Y'|$, there exists $z\in Y'\setminus N_{G'}(S)$ such that $N_{G'}(z)\subseteq X'\setminus S$. As $\delta(G') \ge t$, we have $|X'\setminus S| \ge t$. 
As $|Y'\setminus N_{G'}(S)|=|Y'|-|S|+1 =|X'|-|S|+1\ge t+1$ and $G'$
has at most $t$ vertices of degree less than $\frac{n}{2}$, $Y'\setminus N_{G'}(S)$
contains a vertex of degree at least $\frac{n}{2}$ in $G'$.  As a result, $|X'\setminus S|\ge \frac{n}{2}$,
contradicting the earlier assertion that 
$|X'\setminus S| <\frac{n}{2}$.
 Hence $G'$ has a perfect matching.

There are polynomial time algorithms such as the Hopcroft-Karp algorithm~\cite{matching} in
finding a maximum matching in any bipartite graph, thus a perfect matching of $G$
can be found in polynomial time. 
\qed

\begin{LEM}\label{lem:not-overfull-add-x}
	Let $0<\frac{1}{n_0} \ll  \ve <1$ and $G$ be a star-multigraph with multi-center $x$ on $2n \ge n_0$ vertices.  
	Suppose that $d_G(x)=\delta(G)$, $\delta(G-x) \ge (1+\ve)n$,  and $G-x$ contains no $\Delta(G)$-overfull subgraph.
	Then $G$ contains no $\Delta(G)$-overfull subgraph. 
\end{LEM}

\pf Let  $X\subseteq V(G)$ with $|X|$ odd. We  show that $G[X]$
is not $\Delta(G)$-overfull.  
Assume first that $ 3\le |X| \le 2n-3$. 
By symmetry between $X$ and $V(G)\setminus X$, we further assume $|X| \le n$. 
Then as $\delta(G-x) \ge (1+\ve)n$ and $G-x$ is simple, we have $$e_{G}(X,V(G)\setminus X) \ge (|X|-1)\big((1+\ve)n-|X|+1\big)>2n>\Delta(G).$$
Thus $2e(G[X]) <\Delta(G) |X|-\Delta(G)$ and so $G[X]$
is not $\Delta(G)$-overfull. Thus, the only possible $\Delta(G)$-overfull subgraph of $G$
is obtained from $G$ by deleting a single vertex. 
Note that if $G-u$ is $\Delta(G)$-overfull for some $u\in V(G)$, then so is 
$G-x$ as 
$x$ has minimum degree in $G$. 
However, $G-x$ is not $\Delta(G)$-overfull by our assumption. 
Thus $G-u$ is not $\Delta(G)$-overfull for any $u\in V(G)$. 
Now as $G[X]$ is not $\Delta(G)$-overfull for any $X\subseteq V(G)$ with $|X|$ odd,
it follows that $G$ contains no $\Delta(G)$-overfull subgraph.  
\qed

\subsection{Results on edge colorings}

Let 
$\varphi\in \CC^k(G)$ for  some integer $k\ge 0$. 
For any $v\in V(G)$, the set of colors \emph{present} at $v$ is 
$\varphi(v)=\{\varphi(e)\,:\, \text{$e \in E(G)$ is incident to $v$}\}$, and the set of colors \emph{missing} at $v$ is $\pbar(v)=[1,k]\setminus\varphi(v)$.  For a subset $X$ of $V(G)$ and a color $i\in [1,k]$, define 
$\pbar_X^{-1}(i)= \{v\in X: i\in \pbar(v)\}$. We simply write $\pbar^{-1}(i)$ for $\pbar^{-1}_{V(G)}(i)$. 

 For two distinct colors $\alpha,\beta \in [1,k]$, the components of the subgraph induced by edges with colors $\alpha$ or $\beta$ are called $(\alpha, \beta)$-chains. Each $(\alpha, \beta)$-chain is either a path or an even cycle.   
If we interchange the colors $\alpha$ and $\beta$
on an $(\alpha,\beta)$-chain $C$ of $G$, we get a new edge $k$-coloring  of $G$,  which is denoted by  $\varphi/C$.  This operation is called a \emph{Kempe change}.
For an $(\alpha,\beta)$-chain  $P$, if it is a path with an endvertex $u$, 
we also denote it by $P_u(\alpha,\beta,\varphi)$ to emphasis the endvertex $u$. 
Let  $u,v\in V(G)$.   If $u$ and $u$
are contained in the same  $(\alpha,\beta)$-chain of $G$ with respect to $\varphi$, we say that $u$ 
and $v$ are \emph{$(\alpha,\beta)$-linked} with respect to $\varphi$.

In 1960's, Gupta~\cite{Gupta-67} and, independently, Vizing~\cite{Vizing-2-classes}   provided an upper bound on the chromatic index
of multigraphs, and K\"onig~\cite{MR1511872} gave an exact value 
of the chromatic index for bipartite multigraphs. 

\begin{THM}[\cite{Gupta-67, Vizing-2-classes}]\label{chromatic-index}
Every multigraph  $G$ satisfies $\chi'(G) \le \Delta(G)+\mu(G)$. 
\end{THM}

\begin{THM}[\cite{MR1511872}]\label{konig}
	Every bipartite multigraph $G$ satisfies $\chi'(G)=\Delta(G)$. 
\end{THM}

An edge $k$-coloring of a multigraph $G$ is said to be \emph{equalized} if each color
class contains either $\lfloor |E(G)|/k \rfloor$ or $\lceil |E(G)|/k \rceil$ edges.  McDiarmid~\cite{MR300623} observed the following result. 
\begin{THM}\label{lem:equa-edge-coloring}
	Let $G$ be a multigraph with chromatic index $\chi'(G)$. Then for all $k\ge \chi'(G)$, there is an equalized edge-coloring of $G$ with $k$ colors. 
\end{THM}

Let $G$ be a multigraph, $k\ge 0$ be an integer and $\varphi \in \CC^k(G)$.
There is a polynomial time algorithm to modify $\varphi$ into
an  equalized edge-coloring of $G$ with $k$ colors. To see this, 
suppose $\varphi$
is not equalized and so we take two colors $i,j\in[1,k]$
such that   $\left||\pbar^{-1}(i)|-|\pbar^{-1}(j)|\right|$ 
is largest.  Since $\varphi$
is not equalized,  $\left||\pbar^{-1}(i)|-|\pbar^{-1}(j)|\right| \ge 4$. 
Assume by symmetry that $|\pbar^{-1}(i)|-|\pbar^{-1}(j)| \ge 4$. 
Consider the submultigraph of $G$ induced on the set of edges colored 
by $i$ or $j$, then the submultigraph must have a component that is a path $P$
starting at an edge colored by  $j$ and ending at an edge colored by $j$. 
By swapping the colors $i$ and $j$ along this path $P$, 
we decreased $|\pbar^{-1}(i)|-|\pbar^{-1}(j)|$ by 4. 
Repeating this process, we can obtain an equalized edge-coloring of $G$ with $k$ colors in  $O(k^2 |V(G)|^2)$ time.

\begin{LEM}[\cite{2104.06253}]\label{lem:equi-coloring}
	Let $G$ be a  star-multigraph on $2n$ vertices with multi-center $x$, and let $A$
	and $B$ be a partition of $V(G)$ with $|A|=|B|$, where we assume $x\in A$.  
	If $e(A)=e(B)$, $E_G(A,B) = E_G(x,B)$, and $G$ has an edge coloring  using $k$ colors, 
	then there exists an  edge coloring $\varphi$
	using $k$ colors such that for each $i,j\in [1,k]$, $|\pbar_A^{-1}(i)|=|\pbar_B^{-1}(i)|$
	and $\left||\pbar_A^{-1}(i)|-|\pbar_A^{-1}(j)|\right| \le 2$.  Furthermore,
	such a coloring $\varphi$ can be found in $O(k^2 n^2)$-time.  
\end{LEM}

\begin{LEM}\label{lem:equi-coloring2}
	Let $G$ be a  star-multigraph on $2n$ vertices with multi-center $x$, and let $A$
	and $B$ be a partition of $V(G)$ with $|A|=|B|$, where we assume $x\in A$.  
If $E_G(A,B)=E_G(x,B)$, then  any edge $k$-coloring of $G$ can be modified into 
 an edge $k$-coloring $\varphi$  of $G$  such that the restrictions of $\varphi$ 
 on both $G[A]$ and $G[B]$ are equitable. Furthermore,  such a coloring $\varphi$ can be found in $O(k^2 n^2)$-time.  
  \end{LEM}

\pf  We choose an edge $k$-coloring $\varphi$  of $G$ such that 
$$g_\varphi:=\max_{C\in \{A,B\}}\max_{i,j}\left||\pbar_C^{-1}(i)|-|\pbar_C^{-1}(j)|\right|$$
is smallest and subject to this, the number  $h_\varphi$ of color pairs $(i,j)$
with $\left||\pbar_C^{-1}(i)|-|\pbar_C^{-1}(j)|\right|=g_\varphi$ 
is smallest. 
If $g_\varphi \le 2$, then we are done. Thus $g_\varphi \ge 3$. 
We assume, without loss of generality, that there exist $i,j\in [1,k]$ such that 
$|\pbar_A^{-1}(i)|-|\pbar_A^{-1}(j)|=g_\varphi \ge 3$. 
 As $E_G(A,B) = E_G(x,B)$ and different $(i,j)$-chains are vertex-disjoint,  at most one  vertex from $ \pbar_A^{-1}(i)$ is $(i,j)$-lined with a vertex from $B$. 
As $|\pbar_A^{-1}(i)|-|\pbar_A^{-1}(j)| \ge 3$,  
there exist distinct $u, v \in \pbar_A^{-1}(i)$
such that $P_{u}(i,j,\varphi) = P_{v}(i,j,\varphi)$. 
We now let $\psi=\varphi/P_{u}(i,j,\varphi)$. 
Then $g_{\psi} \le g_\varphi$
and $h_{\psi}<h_{\varphi}$,  
showing a contradiction to the choice of $\varphi$. 
Thus we can find a coloring $\varphi$ such that  $\left||\pbar_C^{-1}(i)|-|\pbar_C^{-1}(j)|\right| \le 2$ for any $C\in \{A,B\}$ and all $i,j\in [1,k]$.

There are  at most ${k \choose 2}$ color pairs,  each color pair takes $O(n)$-time   to reduce the difference of the two corresponding color classes  by 2, 
and $O(n)$ steps to make the two corresponding color classes close in size. Thus a desired edge coloring $\varphi$ can be found in $O(k^2 n^2)$-time. 
\qed 

\begin{LEM}[\cite{2104.06253}]\label{lem:chromatic-index-of-primitive-multiG2}
	If $G$ is a star-multigraph, then $\chi'(G) \le \Delta(G)+1$. 
	Furthermore, an edge coloring of $G$ using at most $\Delta(G)+1$ colors 
	can be found in polynomial time. 
\end{LEM}

A multigraph $G$ is a \emph{near star-multigraph} with multi-center $x$
if, in $G-x$,  there is at most one pair of vertices $u$ and $v$ such that $e_{G-x}(u,v) \ge 2$.

\begin{LEM}\label{lem:chromatic-index-of-primitive-multiG3}
	Let  $G$ be a near star-multigraph with multi-center $x$ such that $\mu(G-x)=e_{G-x}(y,z)$ for some $y,z\in V(G-x)$. Then $\chi'(G) \le  \max\{\Delta(G)+e_G(y,z), \Delta(G)+1\}$. 
	Furthermore, an  edge $\max\{\Delta(G)+e_G(y,z), \Delta(G)+1\}$-coloring of $G$ 
	can be found in polynomial time. 
\end{LEM}

\pf  We may assume  $e_G(y,z) \ge 2$, as otherwise the statement follows from Lemma~\ref{lem:chromatic-index-of-primitive-multiG2}.  
Let $G'$ be obtained from $G$ by deleting $e_G(y,z)-1$ edges joining $y$ and $z$. 
Applying Lemma~\ref{lem:chromatic-index-of-primitive-multiG2} on $G'$, 
we have $\chi'(G') \le \Delta(G')+1 \le \Delta(G)+1$. Now based on any 
edge $(\Delta(G')+1)$-coloring of $G'$, coloring each of the deleted edges using 
a  distinct color gives an edge $(\Delta(G')+1+e_G(y,z)-1)$-coloring of $G$. 
As $\Delta(G')+1+e_G(y,z)-1 \le \Delta(G)+e_G(y,z)$, we have $\chi'(G) \le  \Delta(G)+e_G(y,z)$. 
By Lemma~\ref{lem:chromatic-index-of-primitive-multiG2}, there is a polynomial time algorithm
to edge color $G'$ using at most $\Delta(G)+1$ colors.  Thus there is a polynomial time algorithm
to edge color $G$ using at most $\max\{\Delta(G)+e_G(y,z), \Delta(G)+1\}$ colors. 
\qed 

Given an edge coloring of $G$ and a given color $i$, 
since  vertices presenting $i$
are saturated by the matching consisting of all edges colored by $i$, we have the Parity Lemma below, see~\cite[Lemma 2.1]{MR2028248}.

\begin{LEM}[Parity Lemma]
	Let $G$ be a multigraph and $\varphi\in \CC^k(G)$ for some integer $k\ge \Delta(G)$. 
	Then 
	$|\pbar^{-1}(i)| \equiv |V(G)| \pmod{2}$ for every color $i\in [1,k]$. 
\end{LEM}

\section{Proof of Theorem~\ref{thm:1}}
For $u\in V(G)$, the \emph{deficiency} of $u$
in $G$ is $\df_{G}(u):=\Delta(G)-d_G(u)$. For $U\subseteq V(G)$, 
$\df_{G}(U)=\sum_{u\in U} \df_{G}(u)$.  
We simply write $\df_{G}(V(G))$
as $\df(G)$. For an integer $i \ge 1$,
define $V_i(G)=\{v\in V(G): d_G(v)=i\}$, and we write $V_i$ for $V_i(G)$ if $G$ is clear. 
The proof of Theorem~\ref{thm:1} is based on the following result, which will be proved in the last section.

\begin{THM}\label{thm:D-coloring}
	For all $0<\ve <1$, there exists $n_0$ and $\eta$
	with $0<\frac{1}{n_0} \ll  \eta \ll \ve $
	such that the following statement holds. 
	Let $G$ be a near star-multigraph on $2n\ge n_0$ vertices with multi-center $x$, and let 
	$U=\{v\in V(G)\,:\, \Delta(G)-d_G(v) \ge  \eta n\}$. 
	Suppose that $G$ satisfies  one of the following five conditions: 
		\begin{enumerate}[(a)]
		
		\item $G$ is a star-multigraph and is regular with $ \delta(G) \ge \delta(G-x) \ge (1+\frac{\ve}{2})n $ and $\mu(x)<\eta n$.

		\item $G$ is regular with $\delta(G) \ge (1+\ve )n$, and there exist distinct $y,z\in V(G)\setminus \{x\}$ with the following properties: 
		\begin{enumerate}[(i)]
			\item $d_G^s(x) \ge 2$ and  there are at most $n^{\frac{1}{2}}$ neighbors  $w$ of $x$ satisfying $e_G(x,w) \ge \eta n$;
			\item  $y$ and $z$ are the only two possible vertices in $G-x$
			with $e_{G-x}(y,z) \ge 2$ but  $\mu(G-x) < n^{\frac{1}{2}}$;  and 
			\item $d^s_G(w) \ge (1+\ve)n$ for all $w\in V(G)\setminus\{x\}$.
		\end{enumerate}
			
		\item  $G$ is a star-multigraph and is regular with $\delta(G) \ge (1+\ve )n$, and there exists  $y\in V(G)\setminus \{x\}$ such that 
		\begin{enumerate}[(i)]
			\item $ 2\le d^s(x)< n^{\frac{1}{2}}$; 
			\item  $d^s_G(y) \ge 2\ve n$; and 
			\item $d^s_G(w) \ge (1+\ve)n$ for all $w\in V(G)\setminus\{x, y\}$. 
			
		\end{enumerate}
			\item $G$ is a star-multigraph and  $G$ has two distinct vertices $y,z$ such that 
		\begin{enumerate}[(i)]
			\item $ \mu(x)  <\eta n$; 
			\item  $d_G(y)=d_G(z) =\delta(G)\ge (\frac{1}{2}+\frac{3\ve}{2})n$ and $\min\{d^s_G(y), d_G^s(z)\}>\Delta(G)-\frac{1}{2}(1-0.9\ve)n$; and 
			\item $d_G(w)=\Delta(G) \ge (1+\ve)n$  and $d_{G}^s(w)\ge (1+\ve)n-\mu(x) \ge (1+0.9\ve)n$ for all $w\in V(G)\setminus\{y,z\}$.  
		\end{enumerate}
			\item $G$ is a star-multigraph,  $|U|  \ge \eta n$,  $\delta(G) \ge \delta(G-x) \ge (1+\ve)n$,  and $\mu(x) \le \frac{2}{\eta }$.  
		\end{enumerate}
	Then $\chi'(G)=\Delta(G)$.    
	Furthermore, there is a polynomial time algorithm
	that finds an optimal coloring.     
\end{THM}

\begin{THM1}
	For all $0<\ve <1$, there exists $n_0$
	such that the following statement holds:
	if $G$ is a graph on $2n-1\ge n_0$ vertices with $\delta(G) \ge (1+\ve)n$,  then $\chi'(G)=\Delta(G)$ if and only if $G$ is not overfull. 
	Furthermore, there is a polynomial time algorithm that finds an optimal coloring. 
\end{THM1}

\pf 
The general strategy for proving the first part of the statement  is to add a vertex to $G$ to produce a star-multigraph  $G'$
and then take off some matchings or linear forests from $G'$
such that the resulting multigraph satisfies one of the conditions listed 
in Theorem~\ref{thm:D-coloring}. 
Choose  integer $n_0$ and a real number $\eta $ so that $  0<\frac{1}{n_0} \ll \eta \ll  \ve <1$.

If $\chi'(G)=\Delta(G)$, then clearly $G$ contains no $\Delta(G)$-overfull subgraph. 
Thus we assume that $G$  is not overfull  and show $\chi'(G)=\Delta(G)$. 
Since $n$ is odd and $G$ is not overfull, we have $\df(G) \ge \Delta(G)$.  
We call a  vertex of degree less than $\Delta(G)$ but greater than $\delta(G)$ a \emph{middle degree}  vertex. Let 
$$W=\{v\in V(G)\,:\, \Delta(G)-d_G(v) \ge  \eta n\}.$$ 
We consider the following four  cases regarding the size of $W$.

{\bf \noindent Case 1.  $|W| \ge  2\eta n$. } 

We let $W'\subseteq W$ be a set of $\lfloor \eta n \rfloor $ 
vertices, and let $x$  be a new vertex. 
We form a new graph $G'$ as follows:  add at most $\lfloor \frac{2}{\eta}\rfloor$ edges between $x$ and each vertex of $W'$ one by one until 
the degree of $x$ in the resulting multigraph reaches to $\delta(G)$. 
Since $\lfloor \frac{2}{\eta} \rfloor \times \lfloor  \eta n \rfloor = 2n> \delta(G)$ (we can choose $\eta$ so that $\frac{1}{\eta}$ is an integer), the star-multigraph $G'$ exists.  As $  \eta n>\lfloor\frac{2}{\eta} \rfloor $, we have $\Delta(G')=\Delta(G)$.  By Lemma~\ref{lem:not-overfull-add-x}, 
 $G'$ contains no $\Delta(G')$-overfull subgraph.

For each vertex $u\in V(G')\setminus\{x\}$, we have $d^s_G(u) \ge \delta(G) \ge (1+\ve )n$. 
For each vertex $w\in W\setminus W'$, we have $\Delta(G')-d_{G'}(w)=\Delta(G)-d_G(w) \ge  \eta n$. Furthermore, 
$|W\setminus W'| \ge  \eta n$. 
Now with $W\setminus W'$ playing the role of $U$ in Theorem~\ref{thm:D-coloring}(e) and applying Theorem~\ref{thm:D-coloring}, we see that $\chi'(G')=\Delta(G')=\Delta(G)$ and so $\chi'(G)=\Delta(G)$.

{\bf \noindent Case 2. $|W|=0$.}

The proof of this case is similar to the  beginning part in the proof of Theorem~1.3 in~\cite{2104.06253}, 
but we repeat it to be self-contained.  Note that $|W|=0$ implies that $\Delta(G)-\delta(G)<\eta n$. 
Let $x$ be a new vertex 
and let  $G'$ be obtained from $G$ by adding  $x$ to $G$
and adding some edges between $x$
and vertices $y\in V(G)$ with $d_G(y)<\Delta(G)$
with the following constraints:
\begin{enumerate}[(1)]
	\item $d_{G'}(x)=\delta(G') =\delta(G)$ and $\Delta(G')=\Delta(G)$;   
	\item subject to (1), at most one neighbor of $x$ in $G'$ is not a maximum  degree vertex of $G'$.
\end{enumerate}

Let $V(G)=\{x_1, \ldots, x_{2n-1}\}$ and assume  $d_G(x_1) \le  \ldots \le d_G(x_{2n-1})$.
Note $d_G(x_{2n-1})=\Delta(G)$. 
Let $s\in [1,2n-2]$ be the smallest integer such that $\sum_{i=1}^s\df_G(x_i) \ge \delta(G)$;  such an integer $s$ exists since  $\df(G) \ge \Delta(G)$.  
For each $i\in [1,s-1]$, we add exactly $\df_G(x_i)$ edges between $x$
and $x_i$, and we add $\delta(G)-\sum_{i=1}^{s-1}\df_G(x_i)$ edges between $x$
and $x_s$, and call  $G'$ the resulting multigraph. It is clear that $G'$
satisfies constraints  (1) and (2). Since $\Delta(G)-\delta(G) < \eta n$, 
we have 
$$
|N_{G'}(x)| \ge  \frac{(1+\ve)n}{\eta n} >3. 
$$
By Lemma~\ref{lem:not-overfull-add-x}, 
$G'$ contains no $\Delta(G')$-overfull subgraph.
 
Since $x$ has the smallest degree in $G'$ and $G=G'-x$ is not overfull,  we get 
$\sum_{i\ge 1}\df_{G'}(x_i) \ge \df_{G'}(x)$. Since $|V(G')|=2n$
is even,  $\df_{G'}(x)+\sum_{i\ge 1}\df_{G'}(x_i) $ is even. 
Then by Theorem~\ref{thm:degree-seq}, there exists a multigraph $H$
on $V(G')$ such that $d_H(x)=\df_{G'}(x)$
and $d_H(x_i)=\df_{G'}(x_i)$ for each $i\in [1,2n-1]$.  The multigraph $H$
will aid us in finding a spanning regular subgraph of $G'$. 

Note that $\Delta(H)=\df_{G'}(x)=\Delta(G)-\delta(G) <\eta n$ 
and $H$ contains isolated vertices.
Thus $\chi'(H) \le \Delta(H)+\mu(H) \le 2\Delta(H) < 2\eta n$ by Theorem~\ref{chromatic-index}. Hence we 
can greedily partition $E(H)$ into $k\le 52 \frac{\eta n}{\ve}$ matchings $M_1,\ldots, M_k$ 
each of size at most $\frac{\ve n}{26}$.  Now we take out linear forests (union of vertex-disjoint paths)  
from $G'$ by applying Lemma~\ref{lem:path-decomposition} with $M_1,\ldots, M_k$. 
As at most one neighbor of $x$ is not a maximum degree vertex of $G'$, it follows that 
$x$ is adjacent in $G'$ to at most one vertex of those vertices covered by $M_i$ for each $i\in[1,k]$. 
More precisely, define spanning subgraphs $G_0, \ldots, G_k$
of $G'$ and edge-disjoint linear forests $F_1, \ldots, F_k$
such that 
\begin{enumerate}[(1)]
	\item $G_0:=G'$ and $G_i=G_{i-1}-E(F_i)$ for $i\in [1,k]$,
	\item $F_i$ is a spanning linear forest  in $G_{i-1}$ whose leaves are precisely the 
	vertices in $M_i$. 
\end{enumerate}

Let $G_0=G'$ and suppose that for some $i\in[1,k]$, we already defined $G_0, \ldots, G_{i-1}$
and $F_1, \ldots, F_{i-1}$. As $\Delta(F_1\cup \ldots \cup F_{i-1}) \le 2(i-1) \le 104  \frac{\eta n}{\ve}$,
it follows  that $\delta(G_{i-1}) \ge (1+\ve) n-104  \frac{\eta n}{\ve}\ge (1+\frac{\ve}{2})n$.
Furthermore, by the construction of  $F_1, \ldots, F_{i-1}$, we know that maximum degree vertices of $G'$ are still maximum degree vertices of $G_{i-1}$. Thus in $G_{i-1}$, still 
at most one neighbor of $x$ is not a maximum degree vertex of $G_{i-1}$. As $\mu(x)<\eta n$
and $\delta(G_{i-1})\ge (1+\frac{\ve}{2})n$, $x$ has at least two neighbors in $G_{i-1}$
that are not covered by $M_i$. 
 Since $M_i$ has size at most $\frac{\ve n}{26}$, we can apply Lemma~\ref{lem:path-decomposition} to $G_{i-1}$ and $M_i$ and obtain a spanning linear forest  $F_{i}$ in $G_{i-1}$ whose leaves are precisely the vertices in $M_i$. 
Set $G_i:=G_{i-1}-E(F_{i})$. 

We claim that $G_k$ is regular. Consider any vertex $u\in V(G_k)$. 
For every $i\in[1,k]$, $d_{F_i}(u)=1$ if $u$ is an endvertex of some edge of $M_i$
and $d_{F_i}(u)=2$ otherwise. 
Since $M_1, \ldots, M_k$ partition $E(H)$, we know that $\sum\limits_{i=1}^{k}d_{F_i}(u)=2k-d_H(u)=2k-\df_{G'}(u)$.
Thus 
$$d_{G_k}(u)=d_{G'}(u)-\sum\limits_{i=1}^{k}d_{F_i}(u)=d_{G'}(u)-(2k-\df_{G'}(u))=\Delta(G')-2k.$$
Note that $\mu(x)<\eta n$ and $\delta(G_k-x) \ge (1+\ve)n-104  \frac{\eta n}{\ve}-\mu(x)> (1+\frac{\ve}{2})n$.  
Thus the multigraph $G_k$ satisfies Condition (a) of Theorem~\ref{thm:D-coloring}. 
Now $\chi'(G_k)=\Delta(G_k)$ by Theorem~\ref{thm:D-coloring}. 
We color  the edges of $F_i$ 
using 2 distinct colors from $[\Delta(G')-2k+1, \Delta(G')]$  for each $i\in[1,k]$. 
It is clear that any edge $\Delta(G_k)$-coloring of $G_k$ together with 
this coloring of  $\bigcup_{i=1}^kF_i$ gives an edge coloring 
of $G'$ using $\Delta(G_k)+2k=\Delta(G')$ colors.  Since $G$ is a subgraph of $G'$
and $\Delta(G)=\Delta(G')$, it follows that $\chi'(G)=\Delta(G)$.

%{\bf \noindent Case 3.  $\df(V(G)\setminus W)< \delta(G)$ }. 
%
%Since $\df(G) \ge \Delta(G)$, the condition in Case 3 implies that $W \ne \emptyset$. 
%By Case 1, we have two subcases below to consider. 

{\bf \noindent Case 3. $n^{\frac{1}{2}}\le |W| < 2\eta n$. }  

We add a new vertex $x$. 
For each vertex $v\in V(G)\setminus (V_\Delta \cup W)$, we add $\df_G(v)$ edges between $x$ and $v$ until the degree of $x$ reaches to $\delta(G)$.   If the degree of $x$
reaches to $\delta(G)$ already in this process, we denote the resulting multigraph by $G'$.  
If the degree of $x$ is 
still less than $\delta(G)$ after we adding edges between $x$
and all vertices of $V(G)\setminus (V_\Delta \cup W)$, 
then we add  at most $2n^{\frac{1}{2}}$  edges between $x$
and each vertex of $W$ until the degree of $x$
reaches to $\delta(G)$. Denote the resulting multigraph by $G'$.   As $\df(G) \ge \df_G(W)>\Delta$ and $d_{G'}(x)=\delta(G)$, it follows that $G'$
is not regular.

Note that   $V_\delta(G') \subseteq W \cup\{x\}$  by the construction of $G'$ and so $|V_\delta(G')| \le 2\eta n+1$. 
Thus for any $v\in V(G')\setminus V_\delta(G')$,  $\delta(G'-v-V_\delta(G')) \ge n$
and so both $G'-V_\delta(G')$ and $G'-v-V_\delta(G')$ are hamiltonian by Theorem~\ref{thm:Dirac}. 
Thus if $G'-V_\delta(G')$ and $G'-v-V_\delta(G')$ have even order, then they each has a perfect matching. 
Hence if $|V_\delta(G')|$ is even, we can decrease $\Delta(G')-\delta(G')$ but preserve $\delta(G')$
in deleting a prefect matching  $M$ of $G'-V_\delta(G')$. If $|V_\delta(G')|$ 
is odd but $G'$ has a middle degree vertex $v$, we can decrease $\Delta(G')-\delta(G')$ but preserve $\delta(G')$
in deleting a prefect matching  $M$ of $G'-v-V_\delta(G')$.  Denote by $G_1$
the reduced multigraph from $G'$ by deleting $M$ in either of these two cases. 
We show in the following that we can consider $G_1$ in the place of $G'$. 
As $\delta(G_1-x)=\delta(G'-x) \ge (1+\ve )n$,  $G_1[X]$ is not 
$\Delta(G_1)$-overfull for any $X\subseteq V(G_1)$ with $ 3\le |X|\le 2n-3$.
Thus to show that $G_1$ contains no $\Delta(G_1)$-overfull subgraph, 
we show that $G_1-u$ is not  $\Delta(G_1)$-overfull for any $u\in V(G_1)$. 
If $|V_\delta(G')| \ge 2$,  then for any $u\in V(G_1)$ and $v\in V_\delta(G_1)$, 
we have 
\begin{eqnarray*}
\sum_{w\in V(G_1-u)}\left(\Delta(G_1)-d_{G_1-u}(w)\right)&=&d_{G_1}(u)+(\Delta(G_1)-d_{G_1}(v))+ \df_{G_1}(V(G_1)\setminus \{u,v\})\\
& \ge& \Delta(G_1)+\df_{G_1}(V(G_1)\setminus \{u,v\})\ge \Delta(G_1).
\end{eqnarray*}
Thus $G_1-u$ is not $\Delta(G_1)$-overfull. 
If $|V_\delta(G')|=1$,  then  $V_\delta(G')=\{x\}$.  Note that $x \in V_\delta(G_1)$ and  the middle degree vertex, say $v$,  exists in $G'$. 
Then 
\begin{eqnarray*}
\df(G_1-x) &=&d_{G'}(x)+(\Delta(G_1)-d_{G'}(v))+\df_{G_1}(V(G_1)\setminus \{x,v\})\\
&=&d_{G'}(x)+(\Delta(G')-1-d_{G'}(v))+\df_{G'}(V(G')\setminus \{x,v\})\\
&=&\sum_{w\in V(G'-x)}(\Delta(G')-d_{G'-x}(w))-1. 	
\end{eqnarray*}
 Since $G'$
contains no $\Delta(G')$-overfull subgraph, we have 
$ \sum_{w\in V(G'-x)}(\Delta(G')-d_{G'-x}(w)) \ge \Delta(G')$.  Thus $\df(G_1-x) \ge \Delta(G')-1=\Delta(G_1)$
and so $G_1$ contains no $\Delta(G_1)$-overfull subgraph.
Furthermore, $\chi'(G_1)=\Delta(G_1)$ implies that $\chi'(G')=\Delta(G')$. 
Thus when $|V_\delta(G')|$ is even or $G'$ has a middle degree vertex, we can consider $G_1$ in the place of $G'$ and show that $G_1$
is a class 1 graph.

Thus we assume  that $|V_\delta(G')|$ is odd and $G'$ has no middle degree vertex. This in particular, implies that $\delta(G')$ and $\Delta(G')$ have the same parity. 
As $G'$ has no $\Delta(G')$-overfull subgraph, $|V_\delta(G')|\ge 3$. 
Let $y,z\in V_\delta(G')\setminus\{x\}$ be distinct.  We find a perfect matching $M_{11}$
in $G-(V_\delta(G')\setminus \{y\})$ and a perfect matching $M_{12}$
in $G-(V_\delta(G')\setminus \{z\})$. The matchings exist by Theorem~\ref{thm:Dirac}. 
   Let $G_1=G'-M_{11}-M_{12}$. 
We repeat this same process and find a perfect matching $M_{21}$
in $G_1-(V_\delta(G')\setminus \{y\})$ and a perfect matching $M_{22}$
in $G_1-(V_\delta(G')\setminus \{z\})$. For $i\in [2,(\Delta(G')-\delta(G'))/2]$, we let $G_i=G_{i-1}-M_{i1}-M_{i2}$. Since $d_{G'}(y) \ge \delta(G')+ e_{G'}(x,y) \ge (1+\ve)n+e_{G'}(x,y)$ and $d_{G'}(z) \ge \delta(G')+ e_{G'}(x,z) \ge (1+\ve)n+e_{G'}(x,z)$,  we have 
 $d^s_{G_i}(y), d^s_{G_i}(z) \ge \delta(G')-i$. 
As $\Delta(G')-\delta(G') < 2n-(1+\ve)n=(1-\ve)n$, we see that 
$d^s_{G_i}(y)=d^s_{G_i}(z)=\delta(G')-i \ge (1+\ve)n-\frac{1}{2}(1-\ve)n\ge (\frac{1}{2}+\frac{3\ve}{2})n$. For any vertex $w\in V(G_i)\setminus \{x,y,z\}$, $d^s_{G_i}(w) \ge \delta(G)-\mu(x) \ge (1+0.9\ve)n$. 
As $|V_\delta(G')|+1 <2\eta n+1\ll 0.9\ve n$,  both 
$G_i-(V_\delta(G')\setminus \{y\})$ and $G_i-(V_\delta(G')\setminus \{z\})$ 
satisfy the condition of 
 Lemma~\ref{thm:matching3}, and so each has a 
perfect matching. 
Thus for each $i\in [2,(\Delta(G')-\delta(G'))/2]$, we find matchings $M_{i1}$
and $M_{i2}$ respectively 
from $G_{i-1}-(V_\delta(G')\setminus \{y\})$ and $G_{i-1}-(V_\delta(G')\setminus \{z\})$. 
We claim  that $G^*:=G_{(\Delta(G')-\delta(G'))/2}$ satisfies Condition (d) 
of Theorem~\ref{thm:D-coloring}.    We have 
\begin{enumerate}[(i)]
	\item $\mu(x) <\eta n $, as there might be  edges in $G'$  between $x$ and vertices of $V(G)\setminus (V_\Delta\cup W)$; 
	\item $d_{G^*}(y)=d_{G^*}(z)=\delta(G^*) \ge (\frac{1}{2}+\frac{3\ve}{2})n$; 
	
	\item $\min\{d_{G^*}^s(y), d_{G^*}^s(z)\} \ge  \Delta(G^*)-(\Delta(G')-\delta(G'))/2-\mu(x) \ge \Delta(G^*)-\frac{1}{2}(1-\ve)n-\mu(x)>\Delta(G^*)-\frac{1}{2}(1-0.9\ve)n$; 
		\item  $d_{G^*}(w)=\Delta(G^*) \ge (1+\ve)n$ and  $d_{G^*}^s(w)\ge (1+\ve)n-\mu(x) \ge (1+0.9\ve)n$ for all $w\in V(G^*)\setminus\{y,z\}$. 
\end{enumerate}
By Theorem~\ref{thm:D-coloring}, $\chi'(G^*)=\Delta(G^*)=\delta(G')$. 
Taking an edge $\delta(G')$-coloring of $G^*$,  coloring edges in $M_{i1}$  with color $\delta(G')+2i-1$
and coloring edges in $M_{i2}$  with color $\delta(G')+2i$ for each $i\in [1, (\Delta(G')-\delta(G'))/2]$, 
we obtain an edge $\Delta(G')$-coloring of $G'$ and so an edge $\Delta(G)$-coloring of $G$.

{\bf \noindent Case 4.  $0<|W| < n^{\frac{1}{2}}$. } 

We may assume that $|V_\delta| \ge 2$. Otherwise,  suppose $|V_\delta|=1$; we let 
$M$ 
be a perfect matching of $G-V_\delta$ and let $G_1=G-M$.  
As  $\df(G) \ge \Delta(G)$ 
by our assumption, we have $\df(G_1) = \df(G)-1\ge \Delta(G)-1= \Delta(G_1)$. 
We claim that $G_1$ contains no $\Delta(G_1)$-overfull subgraph. 
As $\delta(G_1) \ge (1+\ve)n$,   $G_1[X]$
 is not $\Delta(G_1)$-overfull for any $X\subseteq V(G_1)$ with $ 3\le |X| \le 2n-3$.
 Thus the only possible $\Delta(G_1)$-overfull subgraph in $G_1$ is $G_1$ itself.   
 However $\df(G_1) \ge \Delta(G_1)$ implies that $G_1$ is not overfull. 
 Since $\delta(G_1)=\delta(G)$,
and $\chi'(G_1)=\Delta(G_1)$ implies that $\chi'(G)=\Delta(G)$, we can consider $G_1$ in the place of $G$. 

Next we may assume that $\df(G) \ge  \Delta(G)+|W|+1$. 
For otherwise, we let $y,z\in V_\delta$  be distinct and add $(\df(G)-\Delta(G))/2$ edges between $y$ and $z$ to get $G^*$ (note that $\df(G)-\Delta(G)$ is even as 
$\df(G)=(2n-1)\Delta(G)-2|E(G)|$), and now add a new vertex $x$ and $\Delta(G)$ edge  between $x$ and the vertices of degree less $\Delta(G)$ in $G^*$
to get a $\Delta$-regular graph $G'$. 
Note that $G'$ may not be star-multiple, but all the possible multiple edges in $G'-x$
are between $y$ and $z$ and $e_{G'}(y,z) \le \frac{1}{2}|W|+1 < n^{\frac{1}{2}}$.  As $\df_{G^*}(w)<\Delta(G^*)$
for any $w\in V(G^*)$, we know that $d_{G'}^s(x) \ge 2$.  By Theorem~\ref{thm:D-coloring}(b), 
we have $\chi'(G')=\Delta(G')=\Delta(G)$ and so $\chi'(G)=\Delta(G)$.

Thus $\df(G)  \ge  \Delta(G)+|W|+1$.  Note that $V_\delta \subseteq W$
and so $|V_\delta| \le |W|<n^{\frac{1}{2}}$. 
We may further assume that $|V_\delta|$ is even and $G$ has no middle degree vertex.
For otherwise, if $|V_\delta|$ is odd, then $G-V_\delta$ has a perfect matching $M$;
if $|V_\delta|$ is even and $G$ has a middle degree vertex $v$, then $G-v-V_\delta$ 
has a perfect matching $M$. In both cases, let $G_1=G-M$.  As $\df(G)  \ge  \Delta(G)+|W|+1$,
we have $\df(G_1)  \ge  \Delta(G)-|W|-1 \ge \Delta(G)>\Delta(G_1)$. Thus $G_1$ is not overfull. 
Furthermore, as $\delta(G_1)=\delta(G) \ge (1+\ve)n$, $G_1[X]$
is not $\Delta(G_1)$-overfull for any $X\subseteq V(G_1)$ with $ 3\le |X| \le 2n-3$. 
Thus $G_1$ contains no $\Delta(G_1)$-overfull subgraph. 
Furthermore $\chi'(G_1)=\Delta(G_1)$ implies that $\chi'(G)=\Delta(G)$. 
Thus we can consider $G_1$ in the place of $G$.

Thus we may assume that  $|V_\delta| \ge 2$ is even, $\df(G)  \ge  \Delta(G)+|W|+1$, and $G$ has no middle degree vertex. 
Let $x$ be a new vertex and $y\in V_\delta$. We first add $\df_G(y)$
edges between $x$ and $y$ and denote the resulting graph by $G'$. 
Note that $V_\delta(G')=\{x\}$. 
  Denote by $\delta'$ the second minimum degree of $G'$. 
Then  we add edges between $x$ and  vertices of $V_{\delta'}(G')=V_\delta \setminus \{y\}$ as follows: 
\begin{enumerate}
	\item [Step 1] If  $|V_{\delta'}(G')|+d_{G'}(x) \le \delta'+1$, then add an edge between $x$ and  each vertex from $V_{\delta'}(G')$ and update $G'$ to be the resulting multigraph; otherwise, $|V_{\delta'}(G')|+d_{G'}(x) \ge \delta'+2$, then we add  an edge
	between $x$ and  each vertex of a selected set of $ \delta'-d_{G'}(x)$ vertices from  $V_{\delta'}(G')$ and update $G'$ to be the resulting multigraph. 
	\item [Step 2]  If $|V_\delta(G')| \ge 2$, then stop; otherwise, return to Step 1. 
\end{enumerate}

Since $\df(G) > \Delta(G)$ and $|V_\Delta| \ge 1$, the process above 
stops before or when the degree of $x$ reaches to $\Delta(G)$. 
Thus $\Delta(G')=\Delta(G)$. Furthermore, we have $d_{G'}(x)=\delta(G') \ge \delta(G)$ and $d^s_{G'}(x) \ge 2$. 
  Lastly by the process,  we have that $x\in V_\delta(G')$,  $|V_\delta(G')| \ge 2$, and $ 0\le d_{G'}(w)- d_{G'}(x) \le 1$ for any $w\in N_{G'}(x)\setminus\{y\}$.
  Actually, as $\df_G(y)+|W| \le (1-\ve)n+n^{\frac{1}{2}}<\delta(G)$, the construction of $G'$
  implies  $N_{G'}(x)=W$.
  As $G=G'-x$ is not overfull and $d_{G'}(x)=\delta(G)$, 
  we know that $G'-u$ is not overfull for any $u\in V(G')$. 
  Also, as $\delta(G'-x) \ge (1+\ve)n$, by calculation, $G'[X]$
  is not $\Delta(G')$-overfull for any $X\subseteq V(G')$ with $ 3\le |X| \le 2n-3$. 
  Thus $G'$ contains no $\Delta(G')$-overfull subgraph. 
  If $G'$
is regular, then  as $ 2\le d^s_{G'}(x) \le |W|$ and  $d_{G'}^s(w) \ge (1+\ve)n$
for any $w\in V(G')\setminus\{x\}$, 
we are done by Theorem~\ref{thm:D-coloring}(c). 

Thus we assume that $G'$ is not regular.  Since $G$ has no middle degree vertex, we have 
 $|V(G')\setminus V_\Delta(G')| \le |(V_\delta(G) \setminus \{y\}) \cup \{x\}| <n^{\frac{1}{2}}$.
Since $\delta^s(G'-x) \ge (1+\ve)n$,
it follows that 
for any $v\in V(G')\setminus V_\delta(G')$,  in the graph  $G'-v-V_\delta(G')$,  every vertex  
 has degree at least $n+1$ except possibly for the vertex $y$, 
 which has degree at least $(1+\ve)n-(2n-(1+\ve)n)-|W|-1>1$. 
Thus both $G'-V_\delta(G')$ and $G'-v-V_\delta(G')$ have a perfect matching  by Lemma~\ref{thm:matching3}. 
Thus if $G'-V_\delta(G')$ and $G'-v-V_\delta(G')$ have even order, then they each has a perfect matching. 
Hence if $|V_\delta(G')|$ is even, we can decrease $\Delta(G')-\delta(G')$ but preserve $\delta(G')$
in deleting a prefect matching  $M$ of $G'-V_\delta(G')$. If $|V_\delta(G')|$ 
is odd but $G'$ has a middle degree vertex $v$, we can decrease $\Delta(G')-\delta(G')$ but preserve $\delta(G')$
in deleting a prefect matching  $M$ of $G'-v-V_\delta(G')$.  Denote by $G_1$
the reduced graph from $G'$ by deleting $M$ in either of these two cases. 
As $\delta(G_1)=\delta(G') \ge (1+\ve )n$, by the same argument as in Case 3,  we know that  $G_1$ contains no $\Delta(G_1)$-overfull subgraph.
Furthermore, $\chi'(G_1)=\Delta(G_1)$ implies that $\chi'(G')=\Delta(G')$. 
Thus we can consider $G_1$ in the place of $G'$ and show that $G_1$
is a class 1 graph.

Thus  $G'$ has no middle degree vertex. This in particular, implies that $\delta(G')$ and $\Delta(G')$ have the same parity.  
Also, as $ 0\le d_{G'}(w)- d_{G'}(x) \le 1$ for any $w\in N_{G'}(x)\setminus \{y\}$, 
it follows that $ d_{G'}(w)=d_{G'}(x) =\delta(G')$ for any $w\in N_{G'}(x)\setminus \{y\}$.
Since  $W=V_\delta = N_{G'}(x)$, we know that   $ d_{G'}(w)=d_{G'}(x) =\delta(G')$ for any $w\in V_\delta\setminus \{y\}$. 
As $V(G)=V_\Delta \cup V_\delta$, and $d_{G'}(y)=\Delta(G')$, it follows that 
$V_\delta(G')=(V_\delta\setminus \{y\}) \cup \{x\}$ and so $|V_\delta(G')|$
is even and $|V_\Delta(G')|$ is even.  
As $\Delta(G')-\delta(G') \le (1-\ve)n$
and $d^s_{G'}(y) \ge (1+\ve)n$,  
by Lemma~\ref{thm:matching3}, $G'-V_\delta(G')$
has a perfect matching $M_1$.  
Let $G_1=G'-V_\delta(G')-M_1$. By Lemma~\ref{thm:matching3} again, 
we can find a perfect matching  $M_2$ in $G_1$. 
In general, we let $G_i=G_{i-1}-M_i$
for each $i\in [2, \Delta(G')-\delta(G')]$. 
We claim that $G_i$ has a perfect matching $M_i$. 
Note that $d_{G_i}(y) \ge (1+\ve )n-(\Delta(G')-\delta(G'))-|V_\delta(G')| \ge (1+\ve)n-(1-\ve)n-n^{\frac{1}{2}} \ge \ve n>1$ and $d_{G_i}(z) \ge \delta(G')-|V_\delta(G')| \ge (1+\frac{\ve}{2})n$. Thus by Lemma~\ref{thm:matching3}, $G_i$ has a perfect matching $M_{i+1}$. 
The graph $G^*:=G_{\Delta(G')-\delta(G')}$ is $\delta(G')$-regular satisfying 
the following conditions:
\begin{enumerate}[(i)]
	\item $2\le d^s_{G^*}(x) \le |W|<n^{\frac{1}{2}}$; 
	\item  $d^s_{G^*}(y) \ge (1+\ve )n-(\Delta(G')-\delta(G')) \ge 2\ve n$; and 
	\item $d^s_{G^*}(w) \ge (1+\ve)n$ for all $w\in V(G^*)\setminus\{x, y\}$, because 
	$w\not\in N_{G'}(x)$ if $w\in V(G^*)\setminus (W\cup \{x\})$
	and $d_{G^*}^s(w)=d_{G'}^s(w) \ge (1+\ve)n$ for any $w\in W\setminus\{y\}$. 
\end{enumerate}
Now by Theorem~\ref{thm:D-coloring}(c), $\chi'(G^*)=\Delta(G^*)$. 
An  edge $\chi'(G^*)$-coloring of $G^*$ together with 
a coloring  with the edges in $M_i$  colored by $i$
for each $i\in [\Delta(G^*)+1,\Delta(G')]$ gives an edge $\Delta(G')$-coloring of $G'$
and so of $G$.

We lastly check that the procedure above yields a polynomial time algorithm. If $G$ is overfull, then $\chi'(G)=\Delta(G)+1$ and 
$G$ can be edge colored using $\Delta(G)+1$ colors in polynomial
time~\cite{MR1156837}.   Thus $G$ is not overfull. 
In each of the Cases 1 to 4, 
it is polynomial time (all the results applied in the proof give appropriate
running time statements) to reduce $G$ 
into a multigraph that satisfies one of the Conditions (a)
to (e) as listed in Theorem~\ref{thm:D-coloring}. 
Then  we find an 
edge $\Delta(G)$-coloring of $G$  in polynomial time by Theorem~\ref{thm:D-coloring}. 
\qed

\section{Proof of Theorem~\ref{thm:D-coloring}}

We will need the following result, which was proved using 
Chernoff bound. 

\begin{LEM}[\cite{2104.06253}, Lemma 3.2]\label{lem:partition}
	There exists a positive integer $n_0$ such that for all $n\ge n_0$ the
	following holds. Let $G$ be a graph on $2n$ vertices, and $N=\{x_1,y_1,\ldots, x_t,y_t\}\subseteq V(G)$, where $1\le t \le n$  is an integer. 
	Then $V(G)$ can be partitioned into two  parts 
	$A$ and $B$ satisfying the properties below:
	\begin{enumerate}[(i)]
		\item  $|A|=|B|$;
		\item $|A\cap \{x_i,y_i\}|=1$ for each $i\in [1,t]$;
		\item $| d_A(v)-d_B(v)| \le n^{\frac{2}{3}}-1$ for each $v\in V(G)$.
	\end{enumerate}
Furthermore, one such partition can be constructed in $O(2n^3 \log_2 (2n^3))$-time. 
\end{LEM}

\begin{THM2}
	For all $0<\ve <1$, there exists $n_0$ and $\eta$
	with $0<\frac{1}{n_0} \ll  \eta \ll \ve $
	such that the following statement holds. 
	Let $G$ be a near star-multigraph on $2n\ge n_0$ vertices with multi-center $x$, and let 
	$U=\{v\in V(G)\,:\, \Delta(G)-d_G(v) \ge  \eta n\}$. 
	Suppose that $G$ satisfies  one of the following five conditions: 
	\begin{enumerate}[(a)]
		
		\item $G$ is a star-multigraph and is regular with $ \delta(G) \ge \delta(G-x) \ge (1+\frac{\ve}{2})n $ and $\mu(x)<\eta n$. 
		\label{Ca}

		\item $G$ is regular with $\delta(G) \ge (1+\ve )n$, and there exist distinct $y,z\in V(G)\setminus \{x\}$ with the following properties: 
		\begin{enumerate}[(i)]
			\item $d_G^s(x) \ge 2$ and  there are at most $n^{\frac{1}{2}}$ neighbors  $w$ of $x$ satisfying $e_G(x,w) \ge \eta n$;
			\item  $y$ and $z$ are the only two possible vertices in $G-x$
			with $e_{G-x}(y,z) \ge 2$ but  $\mu(G-x) < n^{\frac{1}{2}}$;  and 
			\item $d^s_G(w) \ge (1+\ve)n$ for all $w\in V(G)\setminus\{x\}$.
		\end{enumerate}
		\label{Cb}
		
		\item  $G$ is a star-multigraph and is regular with $\delta(G) \ge (1+\ve )n$, and there exists  $y\in V(G)\setminus \{x\}$ such that 
			\begin{enumerate}[(i)]
			\item $ 2\le d^s(x)< n^{\frac{1}{2}}$; 
			\item  $d^s_G(y) \ge 2\ve n$; and 
			\item $d^s_G(w) \ge (1+\ve)n$ for all $w\in V(G)\setminus\{x, y\}$. 
			
		\end{enumerate}
		\label{Cc}
		\item $G$ is a star-multigraph and  $G$ has two distinct vertices $y,z$ such that 
	\begin{enumerate}[(i)]
	\item $ \mu(x)  <\eta n$; 
	\item  $d_G(y)=d_G(z) =\delta(G)\ge (\frac{1}{2}+\frac{3\ve}{2})n$ and $\min\{d^s_G(y), d_G^s(z)\}>\Delta(G)-\frac{1}{2}(1-0.9\ve)n$; and 
	\item $d_G(w)=\Delta(G) \ge (1+\ve)n$  and $d_{G}^s(w)\ge (1+\ve)n-\mu(x) \ge (1+0.9\ve)n$ for all $w\in V(G)\setminus\{y,z\}$.  
\end{enumerate}
		\label{Cd}
		\item $G$ is a star-multigraph,  $|U|  \ge \eta n$,  $\delta(G) \ge \delta(G-x) \ge (1+\ve)n$,  and $\mu(x) \le \frac{2}{\eta }$.  
		\label{Ce}
	\end{enumerate}
	Then $\chi'(G)=\Delta(G)$.    
	Furthermore, there is a polynomial time algorithm
	that finds an optimal coloring.     
\end{THM2}

\pf   We choose $n_0$ and $\eta$ so that 
 $0<\frac{1}{n_0} \ll  \eta \ll \ve $.   
 
 Let  $x_1=x$,  and $y_1\in V(G)\setminus \{x\}$.   Define 
 \begin{numcases}{N^b(x)=}
 \emptyset & \text{Conditions \eqref{Ca}, \eqref{Cd}, and \eqref{Ce}};  \nonumber \\
 \{v\in V(G)\setminus\{y_1\}: e_G(x,v) \ge \eta n\} & \text{Condition \eqref{Cb}}; \nonumber \\
 N_G(x)\setminus\{y_1\} &  \text{Condition \eqref{Cc}}.\nonumber
 \end{numcases}
  We then define a set $N$ according to the different conditions as follows:  
 \begin{enumerate}[]
 		\item Condition~\eqref{Ca}: let  $N=\{x_1,y_1\}$.  
 		\item Condition~\eqref{Cb}: let $N^b(x)=\{x_2, \ldots, x_t\}$.   By the assumption of Condition \eqref{Cb}, we have  $t-1<n^{\frac{1}{2}}$. 
 		We take distinct $y_2,\ldots, y_t \in V(G)\setminus \{x_1, y_1, x_2,\ldots, x_t\}$ 
 		and define $N=\{x_1,y_1, \ldots, x_t, y_t\}$.    
 	\item Condition~\eqref{Cc}: let $N_G(x)\setminus\{y_1\}=\{x_2, \ldots, x_t\}$. By the assumption of Condition \eqref{Cc}, we have  $t-1<n^{\frac{1}{2}}$. We take distinct $y_2,\ldots, y_t \in V(G)\setminus \{x_1, y_1, x_2,\ldots, x_t\}$ 
 	and define $N=\{x_1,y_1, \ldots, x_t, y_t\}$.

 	\item  Condition~\eqref{Cd}: let  $N=\{x_1,y_1,x_2,y_2\}$,  where $x_2=y$
 	and $y_2=z$.  
 	
 		\item 	Condition~\eqref{Ce}: take $2\lfloor(2n-|V_\Delta|)/2 \rfloor$ vertices from $V(G)\setminus V_\Delta$ 
 	and name them as $x_2,y_2, \ldots, x_t, y_t$, where $t-1:=\lfloor (2n-|V_\Delta|)/2 \rfloor$ and  we assume that the first $\lfloor |U|/2\rfloor$ pairs of vertices $x_i,y_i$ are all from $U$.
 	Let $N=\{x_1,y_1, \ldots, x_t, y_t\}$. 
 \end{enumerate}
We say that $y_i$ and $x_i$ from the set $N$ are \emph{partners} of each other.

Applying Lemma~\ref{lem:partition} on the underlying  graph of $G$ and $N$,
we obtain a partition $\{A, B\}$ of $V(G)$ satisfying the following properties in the underlying graph of $G$:  $|A|=|B|$;  $|A\cap \{x_i,y_i\}|=1$ for each $i\in [1,t]$; and 
$| d_{A}(v)-d_{B}(v)| \le n^{\frac{2}{3}}-1$ for each $v\in V(G)$.

For the pair $(x_1,y_1)$, where recall $x_1=x$, 
by switching the two vertices from their current partitions  and  renaming the resulting partitions if necessary,  we may assume that   $x\in A$ and $d_B(x) \ge d_A(x)$. 
When $G$ is in Conditions~\eqref{Cb} or~\eqref{Cc}, if $N^b(x)\cap A \ne \emptyset$,  
we move vertices in $N^b(x)\cap A$ from $A$ to $B$ and move  the partners of vertices 
in $N^b(x)\cap A$ from $B$ to $A$.  
Still denote by $A$ and $B$ the 
resulting partition of $V(G)$.  Then we still have $d_B(x) \ge d_A(x)$
as $e_G(x,x_i) \ge e_G(x,y_i)$ for any $x_i\in N^b(x)$.  
Considering now  the partition $\{A,B\}$ in the multigraph $G$, it has the properties below, 
where ``Pa'' indicates the property when $G$ is in Condition~\eqref{Ca}, and similar meaning for 
``Pb'' to  ``Pe.'' 
\begin{enumerate}[Pa]
	\item  \begin{eqnarray}
		|d_{A}(v)-d_{B}(v)| &\le& n^{\frac{2}{3}}+ e_G(x,v) \quad \text{for each $v\in V(G)\setminus\{x\}$.} \label{partition-Ca} \nonumber 
	\end{eqnarray}
		\label{Pa} 
	\item \begin{eqnarray}
		N^b(x) &\subseteq& B; \nonumber \\ 
		|d_{A}(v)-d_{B}(v)| &\le& n^{\frac{2}{3}}+2n^{\frac{1}{2}}+n^{\frac{1}{2}}+e_G(x,v) \quad \text{for each $v\in V(G)\setminus\{x\}$,} \label{partition-Cb}  \nonumber
	\end{eqnarray}
where the term ``$n^{\frac{1}{2}}$'' is to take account the the possible number of multiple edges 
joining $x$ and $z$. 
\label{Pb}

\item \begin{eqnarray}
	N_G(x) &\subseteq& B; \nonumber \\ 
	|d_{A}(v)-d_{B}(v)| &\le& n^{\frac{2}{3}}+2n^{\frac{1}{2}}+e_G(x,v) \quad \text{for each $v\in V(G)\setminus\{x\}$.} \label{partition-Cc}  \nonumber
\end{eqnarray}
\label{Pc}
\item  \begin{eqnarray}
	|\{y,z\}\cap A| &=& |\{y,z\}\cap B|=1; \nonumber \\ 
	|d_{A}(v)-d_{B}(v)| &\le& n^{\frac{2}{3}}+e_G(x,v) \quad \text{for each $v\in V(G)\setminus\{x\}$.} \label{partition-Cd} \nonumber
\end{eqnarray}
\label{Pd}

\item  \begin{eqnarray}
	|U\cap A| \ge \lfloor \frac{|U|}{2} \rfloor,  & & |U\cap B| \ge \lfloor \frac{|U|}{2} \rfloor; \nonumber \\ 
	|d_{A}(v)-d_{B}(v)| &\le& n^{\frac{2}{3}}+ e_G(x,v) \quad \text{for each $v\in V(G)\setminus\{x\}$.} \label{partition-Ce} \nonumber
\end{eqnarray}
\label{Pe}
\end{enumerate}

Let 
$$ G_A=G[A], \quad G_B=G[B], \quad \text{and } \quad H=G[A,B]. $$

Define $G_{A,B}$ to be the union of  $G[A]$, $G[B]$ together with $ \lfloor (d_B(x)-d_A(x))/2 \rfloor$ edges incident with $x$
from $E(H)$.  
 When $G$ is in Conditions~\eqref{Cb} or~\eqref{Cc},   we take at most  $\lceil e_G(x,u)/2\rceil $ edges
between $x$ and each $u\in B$.   In particular, when  $G$ is in Conditions~\eqref{Cc}, 
as $N_G(x) \subseteq B$, we take either $\lfloor e_G(x,u)/2 \rfloor$   or $\lceil e_G(x,u)/2\rceil$ edges between $x$ and each $u\in N_G(x)$.   
By this definition of $G_{A,B}$ and the Properties~P\ref{Pa} to~P\ref{Pe} above,  we have 
\begin{eqnarray} 
	d_{G_{A,B}}(x) &=&  \lfloor d_G(x)/2 \rfloor,  \nonumber \\
d_{G_{A,B}}(v) &\ge& \frac{1}{2}\left(d_G(v)-e_G(v,x)-n^{\frac{2}{3}}-2n^{\frac{1}{2}}-n^{\frac{1}{2}}\right)+e_G(v,x) \quad  \text{for 	$v\in A\setminus\{x\}$,}  \label{eqn:degree-v-in-Gab} \\
d_{G_{A,B}}(v) &\ge& \frac{1}{2}\left(d_G(v)-e_G(v,x)-n^{\frac{2}{3}}-2n^{\frac{1}{2}}-n^{\frac{1}{2}}\right)\quad  \text{for 	$v\in B$.}  \nonumber 
\end{eqnarray}
For vertex $v\in V(G)\setminus(N_b(x)\cup \{x\})$, as $e_G(x,v)<\eta n$ when $G$ is in 
Conditions~\eqref{Ca} to~\eqref{Cd} and $e_G(x,v) \le \frac{2}{\eta}$ when $G$ is in 
Condition~\eqref{Ce}, we will use the 
following bound: 
\begin{numcases}{d_{G_{A,B}}(v) \ge}
\frac{1}{2}\left(d_G(v)-1.2\eta n\right) & \text{when $G$ is in 
	Conditions~\eqref{Ca} to~\eqref{Cd};}  \nonumber \\ 
\\ 
\frac{1}{2}\left(d_G(v)-2n^{\frac{2}{3}}\right) & \text{when $G$ is in 
	Condition~\eqref{Ce}.} \nonumber 
\end{numcases}

For $v\in N^b(x)$, we have 
\begin{equation}\label{eqn:general-degree-v-in-Gab2}
	d_{G_{A,B}}(v) \ge \frac{1}{2}\left(d^s_G(v)-2n^{\frac{2}{3}}\right). 
\end{equation}

Furthermore,  when  $G$ is  in Conditions~\eqref{Ca} to~\eqref{Cd},  
by the assumptions on $G$ and the partition $\{A,B\}$, 
 we have $\sum_{v\in A}d_G(v)=\sum_{v\in B}d_G(v)$ and so $e(G_A)=e(G_B)$.

To prove the theorem, we will construct an edge coloring 
of $G$ using $\Delta(G)$ colors. 
We provide below an overview of the steps. 
At the start of the process, $E(G)$ is assumed to be uncolored,
and throughout the process, the partial edge coloring of $G$
is always denoted by $\varphi$, which is updated step by step. 

\begin{enumerate}[Step 1]
	\item Define $S=\{v\in V(G)\,:\, \Delta(G)-d_G(v) \ge 7n^{\frac{2}{3}}\}$. 
	Let $k=\Delta(G_{A,B})+ \lfloor  n^{\frac{1}{2}} \rfloor$. 
	Since  $G$ is a near star-multigraph with $\mu(G-x)<n^{\frac{1}{2}}$, 
	by Theorem~\ref{lem:chromatic-index-of-primitive-multiG3}, we find an edge $k$-coloring  $\varphi$ of  $G_{A,B}$. 
	If there exist distinct $u,v\in S\cap A$  or distinct $u,v \in S\cap B$ such that $\pbar(u)\cap \pbar(v)\ne \emptyset$, we add an edge joining $u$
	and $v$ and color the new edge by a color in $\pbar(u)\cap \pbar(v)$. The edge coloring $\varphi$ is updated and we still call it $\varphi$. 
	We iterate this process of adding and coloring edges and call the  multigraphs  resulting from $G_A, G_B$ and  $G_{A,B}$, respectively, as $G_A^*$, $G_B^*$, and $G^*_{A,B}$, and call 
	$G^*$ the union of $G^*_{A,B}$ and $H$. It can be verified that $\Delta(G^*)=\Delta(G)$.  We will modify the current 
	edge coloring, which is still named $\varphi$, such that the following properties are satisfied:
	\begin{enumerate}
		\item[S1.1 ]  Let $i,j\in [1,k]$ and $C\in \{A,B\}$. When $G$ is in Conditions~\eqref{Ca} to~\eqref{Cd}, we have  $$
		|\pbar^{-1}_A(i)|=|\pbar^{-1}_B(i)| \quad \text{and} \quad  ||\pbar^{-1}_C(i)|-|\pbar^{-1}_C(j)|| \le 2\quad 
		$$
		When $G$ is in Condition~\eqref{Ce}, we have 
			$$	\left||\pbar_C^{-1}(i)|-|\pbar_C^{-1}(j)|\right| \le 2.$$
			\label{s1.1}
		\item[S1.2 ] Let  $i\in [1,k]$.  When $G$ is in Conditions~\eqref{Ca} to~\eqref{Cd}, we have 
		$$
		|\pbar^{-1}_A(i)| < 3 \eta n \quad \text{and} \quad  |\pbar^{-1}_B(i)| < 3\eta n.
				$$
				When $G$ is in Condition~\eqref{Ce}, we have 
		$$
		|\pbar^{-1}_A(i)| < 7 n^{\frac{2}{3}} \quad \text{and} \quad  |\pbar^{-1}_B(i)| < 7 n^{\frac{2}{3}}.
		$$
		\label{s1.2}
	\end{enumerate}

	\item Modify the partial edge-coloring of $G^*$ obtained in Step 1 by exchanging alternating paths. Upon the completion of Step 2, each of the $k$ color class will be a 1-factor of $G^*$. In the process of Step 2,  a few edges of $H$ will be colored and  a few edges of $G^*_A \cup G^*_B$ will be uncolored. Denote by  $R_A$ and $R_B$, respectively,  the submultigraphs of $G^*_A$
	and $G^*_B$ consisting of the uncolored edges. The two multigraphs $R_A$ and $R_B$  will initially be empty, but 
	one, two, three, or four edges will be added to at least one of them 
	when each time we exchange an alternating path. The calculations involved when $G$
	is in Condition~\eqref{Ce} is essentially different than the other cases, 
	so we need define parameters  correspondingly as follows: let 
		\begin{numcases}{s=}
		3\eta n^2 & \text{when $G$ is in Conditions~\eqref{Ca} to~\eqref{Cd}}, \nonumber  \\
		7n^{\frac{5}{3}} & \text{when $G$ is in Condition~\eqref{Ce}},  \nonumber 
	\end{numcases}
	and 
\begin{numcases}{r=}
	\eta^{\frac{1}{2}} n & \text{when $G$ is in Conditions~\eqref{Ca} to~\eqref{Cd}}, \nonumber  \\
	n^{\frac{5}{6}} & \text{when $G$ is in Condition~\eqref{Ce}}.   \nonumber 
\end{numcases}

	The conditions below will be satisfied at the completion of this step:
	\begin{enumerate}
		\item [S2.1] The number of uncolored edges in each of $G_A^*$
		and $G_B^*$ is less than $4s$. 
		When $G$ is in Conditions~\eqref{Ca} to~\eqref{Cd}, $G_A^*$
		and $G_B^*$ have the same number of uncolored edges; and when $G$
		is in Condition~\eqref{Ce}, the number of uncolored edges in $G_B^*$
		is greater than or equal to the number of uncolored edges in $G_A^*$ (this follows if we assume by symmetry that  $e(G_A^*) \le e(G_B^*)$). 
		\label{s21}
		\item [S2.2] $\Delta(R_A)$ and $\Delta(R_B)$ are less than $r$.
		\label{s22}
		\item [S2.3]  For $C\in \{A,B\}$,   let $
		S_C=\{u\in S\cap C\,:\, d_{G_{A,B}^*}(u) \le k-2n^{\frac{2}{3}} \}. 
		$
		We require
 every vertex  $u\in V(G^*)\setminus (S_A\cup S_B\cup N^b(x)\cup \{x\})$ is incident in $G^*$ with fewer than $|\pbar(u)|+r$ colored edges of $H$; and every vertex  $u\in  S_A\cup S_B\cup N^b(x)\cup \{x\}$ is incident in $G^*$ with at most $|\pbar(u)|$ colored edges of $H$.

\label{s23}
	\end{enumerate}
	
	\item We will edge color each of $R_A$
	and $R_B$ and color a few uncolored edges of $H$ using another 
	$\ell$ colors, where $\ell =2\lfloor  r   \rfloor$.  
	The goal is to ensure that 
	each of these  $\ell$ new color classes obtained at the completion of Step 3 presents at all vertices from $V(G^*)\setminus U$ while preserve the $k$ 1-factors already obtained through Steps 1 and 2.  
	\item At the start of Step 4, all of the uncolored edges of $G^*$ belong to $H$. Denote by $R$ the subgraph of $G^*$ consisting of the uncolored edges. 
	We will show that $\Delta(R)=\Delta(G^*)-k-\ell$. 
	 This submultigraph is bipartite, so we can color its edges using $\Delta(G^*)-k-\ell$
	colors by Theorem~\ref{konig}.  
\end{enumerate}

When Step 4 is completed, we obtain an edge coloring of $G^*$ using exactly $\Delta(G^*)$ colors. We now give the details of each step,
and for notation that was already defined in the outline above,
we will use them directly.

\bigskip 

\bigskip 

\begin{center}
 	Step 1: Coloring $G_{A,B}$ 
\end{center}

Recall $S=\{v\in V(G)\,:\, \Delta(G)-d_G(v) \ge 7n^{\frac{2}{3}}\}$.
 Note that when $G$ is in Conditions~\eqref{Ca} to~\eqref{Cc}, $S=\emptyset$;
 when $G$ is in Condition~\eqref{Cd}, $S\subseteq\{y,z\}$; 
 and when $G$ is in Condition~\eqref{Ce}, $U\subseteq S$. 
 When $G$ is in Condition~\eqref{Cd}, we have $|\{y,z\}\cap A|=|\{y,z\}\cap B|=1$
 by Property P\ref{Pd}. Thus,  when $G$ is in Conditions~\eqref{Ca} to~\eqref{Cd},
 we have $G_A^*=G_A$, $G_B^*=G_B$, $G_{A,B}^*=G_{A,B}$, and $G^*=G$.

Following the operations described in the outline of Step 1, 
for the current edge coloring $\varphi$ of $G_{A,B}^* $,  
the following statement holds:  $\pbar(u)\cap \pbar(v)=\emptyset$
for any two distinct $u,v\in S\cap A$ or any two distinct $u,v\in S\cap B$.  
Therefore, for any $C\in \{A,B\}$, we have  
 \begin{equation}\label{eqn:sum-of-def-in-S}
 \sum_{u\in C\cap S}|\pbar(u)|=\sum_{u\in C\cap S} (k-d_{G^*_{A,B}}(u)) \le k. 
 \end{equation}
 
 Note that $k \le \Delta(G_{A,B})+n^{\frac{1}{2}}$, and
 by the partition $\{A,B\}$ of $V(G)$, we have 
 $\Delta(G_{A,B}) \le \max_{v\in V(G)} \{\frac{1}{2}\left(\Delta(G)-e_G(x,v)+n^\frac{2}{3}+2n^\frac{1}{2}+ n^{\frac{1}{2}}\right)+e^*\}$, where $e^*=e_G(x,v)<\eta n$ if $v\in A$
 and $e^*= \lceil e_G(x,v)/2 \rceil $  if $v\in B$. 
Thus when $G$ is in Conditions~\eqref{Ca} to~\eqref{Cd},   
 we have 
 \begin{equation}\label{eqn:k-value-a-d}
 	k \le \frac{1}{2}\left(\Delta(G)+n^\frac{2}{3}+2n^\frac{1}{2}+ n^{\frac{1}{2}}\right) +\frac{1}{2}\eta n+n^{\frac{1}{2}}<\frac{1}{2}\Delta(G)+0.6\eta n. 
 \end{equation}
 When $G$ is in Condition~\eqref{Ce}, 
 as $\mu(G)<\frac{2}{\eta}$,  we have 
 \begin{equation}\label{eqn:k-value-e}
 	k \le \frac{1}{2}\left(\Delta(G)+n^\frac{2}{3}+2n^\frac{1}{2}+ n^{\frac{1}{2}}\right) +\frac{1}{\eta}+n^{\frac{1}{2}}<\frac{1}{2}\Delta(G)+n^\frac{2}{3}. 
 \end{equation}
 
We will in the rest of the proof show that $\chi'(G^*)=\Delta(G^*)$,
this is because $G$ is a subgraph of $G^*$ and  $\Delta(G^*)=\Delta(G)$.
To see  $\Delta(G^*)=\Delta(G)$, 
note that the edge additions in Step 1 happen only when $G$ is in Condition~\eqref{Ce}. 
Under Condition~\eqref{Ce},  
for any $u\in S\cap A$, we have 
$$
d_{G^*}(u) \le 
k +e_G(u, B)  \le  \frac{1}{2}\Delta(G)+n^\frac{2}{3}+ +\frac{1}{2}(\Delta(G)-7n^{\frac{2}{3}}+n^{\frac{2}{3}}+\frac{2}{\eta}) < \Delta(G).
$$
Similarly, we have $d_{G^*}(u) \le \Delta(G)$ for any $u\in S\cap B$.  
In particular, if $u\in U$, as $\Delta(G)-d_G(u) \ge \eta n$, 
we have 
\begin{equation}\label{eqn:degree-U-vertex-in-G*}
d_{G^*}(u) \le  \frac{1}{2}\Delta(G)+n^\frac{2}{3}+\frac{1}{2}(\Delta(G)-\eta n+n^{\frac{2}{3}}+\frac{2}{\eta})< \Delta(G)-\frac{\eta n}{3}.
\end{equation}

When $G$ is in Conditions~\eqref{Ca} to~\eqref{Cd},  since $e(G_A^*)=e(G_A)=e(G_B)=e(G_B^*)$, 
by Lemma~\ref{lem:equi-coloring},  we modify 
$\varphi$  into an edge $k$-coloring of $G^*_{A,B}$, and still call  $\varphi$ the resulting 
coloring, such that 
$$
|\pbar^{-1}_A(i)|=|\pbar^{-1}_B(i)| \quad \text{and} \quad  ||\pbar^{-1}_C(i)|-|\pbar^{-1}_C(j)|| \le 2\quad \text{for $i,j\in [1,k]$ and $C\in \{A,B\}$}. 
$$
When $G$ is in Condition~\eqref{Ce}, by Lemma~\ref{lem:equi-coloring2}, we modify $\varphi$ so that 
$$
||\pbar^{-1}_C(i)|-|\pbar^{-1}_C(j)|| \le 2\quad \text{for $i,j\in [1,k]$ and $C\in \{A,B\}$}.
$$ 

Note that under the modified  coloring $\varphi$, it is possible that $\pbar(u)\cap \pbar(v) \ne \emptyset$
for some distinct  $ u,v\in S\cap A$ or distinct $u,v\in S\cap B $. 
However  the inequalities in~\eqref{eqn:sum-of-def-in-S} still hold.  
For $C\in \{A,B\}$, define 
\begin{eqnarray*}
S_C=\{u\in S\cap C\,:\, d_{G_{A,B}^*}(u) \le k-2n^{\frac{2}{3}} \}. 
\end{eqnarray*}
 Next, we verify that every color $i\in [1,k]$ is missing at a small number of vertices. 
We first provide an upper bound on $|\pbar(u)|$ for each $u\in V(G^*)$. 
By the definition of $G_{A,B}$, \eqref{eqn:general-degree-v-in-Gab}, \eqref{eqn:general-degree-v-in-Gab2}, \eqref{eqn:k-value-a-d}, and~\eqref{eqn:k-value-e}, as $|\pbar(u)| = k-d_{G^*_{A,B}}(u)$, we have
   	   	\begin{numcases}{|\pbar(u)| \le}
 		\frac{1}{2}\Delta(G)+0.6\eta n-\lfloor\Delta(G)/2\rfloor< \eta n  
 		& \text{$u=x$, Conditions~\eqref{Ca} to~\eqref{Cd}};\nonumber \\
 		& \nonumber \\ 
 		\frac{1}{2}\Delta(G)+n^{\frac{2}{3}}-\frac{1}{2} (1+\ve)n  < \frac{1}{2}(1-\ve)n+n^{\frac{2}{3}} 
 		& \text{$u=x$, Condition~\eqref{Ce}};  \nonumber \\
 		& \nonumber \\ 
 		\frac{1}{2}\Delta(G)+0.6\eta n- \frac{1}{2}\left(\Delta(G)-7n^{\frac{2}{3}}-1.2\eta n \right) & \text{$u\in V(G)\setminus(S\cup N^b(x)\cup \{x\})$, }\nonumber \\
 		<2\eta n &  \text{and Conditions~\eqref{Ca} to~\eqref{Cd}};\nonumber \\
 		& \nonumber  \\ 
 		\frac{1}{2}\Delta(G)+n^{\frac{2}{3}}-\frac{1}{2}\left(\Delta(G)-7n^{\frac{2}{3}}-
 		2n^{\frac{2}{3}}\right)   &  \text{$u\in V(G)\setminus S$, Condition~\eqref{Ce}}; \nonumber  \\ 
 		<6n^{\frac{2}{3}} & \nonumber \\
 			& \nonumber \\ 
 		 \frac{1}{2}\Delta(G)+0.6\eta n-\left(\frac{1}{2}(1+\ve)n-2n^{\frac{2}{3}}\right) & 
 		\text{$u\in N^b(x)$, and}  \nonumber  \\ 
 		< (\frac{1}{2}-\frac{1}{3}\ve)n & \text{Conditions~\eqref{Cb} and~\eqref{Cc}}; \nonumber \\
 		& \nonumber \\  
 		\frac{1}{2}\Delta(G)+0.6\eta -&  \text{$u=y$, Condition~\eqref{Cc}};\nonumber \\ 
 		 \left(\frac{1}{2}\left(\Delta(G)-e_G(x,u)-2n^{\frac{2}{3}} \right)+\lfloor e_G(x,u)/2 \rfloor \right)  &  \nonumber  \\ 
 		< \eta n &  \nonumber \\
 	 		\frac{1}{2}\Delta(G)+0.6\eta-\frac{1}{2}\left(d^s_G(u)-1.2\eta n \right) & \text{$u\in \{y,z\}$, Condition~\eqref{Cd}}; \nonumber  \\ 
 		=\frac{1}{2}\Delta(G)+1.2\eta-\frac{1}{2}\left(\Delta(G)-\frac{1}{2}(1-0.9\ve)n \right) &  \nonumber \\ 
 		<(\frac{1}{4}-\frac{1}{5}\ve) n  & \nonumber \\
 		& \nonumber \\  
 		\frac{1}{2}\Delta(G)+n^{\frac{2}{3}}-\frac{1}{2}\left((1+\ve)n-2n^{\frac{2}{3}} \right) & \text{$u\in S_A\cup S_B$, Condition~\eqref{Ce}}; \nonumber  \\ 
 		< (\frac{1}{2}-\frac{1}{3}\ve)n&  \nonumber \\
 		& \nonumber \\  
 		2n^{\frac{2}{3}} & \text{$u\in S\setminus(S_A\cup S_B)$ and} \nonumber \\ &\text{Condition~\eqref{Ce}.}\nonumber 
 	\end{numcases}

  Since $|\pbar(u)|<2n$ for each $u\in N^b(x)$ and $|N^b(x)|<n^{\frac{1}{2}}$, 
the inequalities above together  with the fact in~\eqref{eqn:sum-of-def-in-S} give  the following inequalities for any $C\in \{A,B\}$:  
   \begin{eqnarray} \label{eqn:total-missing-color-a-d}
 	\sum\limits_{u\in C}|\pbar(u)| &\le &2nn^{\frac{1}{2}}+ 2\eta n^2+ k < 3\eta n^2-2n \quad \text{$G$ is in Conditions~\eqref{Ca} to~\eqref{Cd};}\\
 	\sum\limits_{u\in C}|\pbar(u)| &\le &2nn^{\frac{1}{2}}+ 6n^{\frac{5}{3}} + k < 7n^{\frac{5}{3}}-2n  \quad \text{$G$ is in Condition~\eqref{Ce}.} \label{eqn:total-missing-color-e}
 \end{eqnarray}

By  SI.1,  
the Statement S1.2 stated in the outline of Step 1 holds.

\begin{center}
	Step 2: Extending  existing color classes into  1-factors  
\end{center}

We will extend each of the $k$ color classes obtained in Step 1  into 
a 1-factor of $G^*$
by exchanging  alternating paths.  Each alternating path  starts and also ends 
at an uncolored edge of $H$
and is alternating between 
 colored edges  of $G_A^*\cup G_B^*$
and uncolored edges of $H$.  Thus 
during the procedure of Step 2, we will uncolor some of the edges of 
$G^*_A\cup G_B^*$, and will color some of the edges of $H$. 
Recall that $R_A$ and $R_B$ are the submultigraphs of $G^*_A$
and $G^*_B$ consisting of the uncolored edges, which are 
 initially empty. 

To ensure Condition S2.2 is satisfied, we say that an edge  $e=uv\in E(G^*_A\cup G_B^*)$ is \emph{good} if $e\not \in E(R_A \cup R_B)$  and  the degree of $u$ and $v$ 
in both $R_A$ and $R_B$ is less than $r$ (actually, note that when  $uv\in E(G_A^*)$, then the degree of $u$
and $v$ is zero in $R_B$  and vice versa). 
Thus  a good edge can be added to $R_A$ or $R_B$ without violating S2.2.

 By S1.\ref{s1.1},  for each color 
 $i\in [1,k]$, we pair up a vertex from $\pbar^{-1}_A(i)$ 
 and a vertex from $\pbar^{-1}_B(i)$,  and then pair up 
 the remaining unpaired vertices from  $\pbar^{-1}_A(i)$ or $\pbar^{-1}_B(i)$  
 as $|\pbar^{-1}_A(i)|-|\pbar^{-1}_B(i)|$  is even by the Parity Lemma. Each of those pairs is called 
 a \emph{missing-common-color pair} or \emph{MCC-pair} in short with respect to the color $i$. 
 In particular, 
 when $G$ is in Conditions~\eqref{Ca} to~\eqref{Cd},  every vertex from $\pbar^{-1}_A(i)$
 is paired up with a vertex from $\pbar^{-1}_B(i)$. 
 
 For every MCC-pair $(a,b)$ with respect to some color $i\in[1,k]$, 
 we will exchange an alternating path $P$ from $a$ to $b$ with at most 13
 edges, where, if exist, the first, third, fifth, seventh,  ninth,  eleventh, and thirteenth edges are uncolored
 and the second, fourth, sixth, eighth,  tenth, and twelfth  edges are good edges colored by $i$. After $P$ is exchanged, $a$ and $b$ will be incident with edges  colored by $i$, and at most four good edges will be added to each of $R_A$ and $R_B$. With this information at hand, 
 before demonstrating the existence of such paths,  we  show that Conditions S2.1, S2.2 and S2.3  can be guaranteed  at the end of Step 2.

 After the completion of Step 1, by~\eqref{eqn:total-missing-color-a-d} and~\eqref{eqn:total-missing-color-e}, the total number 
 of missing colors from vertices in $A$
or from vertices in $B$ is at most $s$, which is $3\eta n^2$
when $G$ is in Conditions~\eqref{Ca} to~\eqref{Cd}, and is $7n^{\frac{5}{3}}$
when   $G$ is in Condition~\eqref{Ce}. 
 Thus there are at most $s$ MCC-pairs. 
 For each MCC-pair $(a,b)$ with $a,b\in V(G^*)$, 
 at most four edges will be added to each of $R_A$ and $R_B$
 when we exchange an alternating path from $a$ to $b$.  Thus 
 there will always be fewer than  
 $4s$
 edges in each of $R_A$ and $R_B$. At the completion of 
 Step 2, each of the $k$ color classes is a 1-factor of $G^*$
 so each of $G_A^*$ and $G_B^*$ have the same number of 
 colored edges. When $G$ is in Conditions~\eqref{Ca} to~\eqref{Cd}, since  $e(G_A^*)=e(G_A)=e(G_B)=e(G_B^*)$,  
 we have $e(R_A)=e(R_B)$. When $G$ is in Condition~\eqref{Ce}, 
 suppose, without loss of generality, that $e(G_B^*) \ge e(G_A^*)$,
 then we get $e(R_B) \ge e(R_A)$. 
  Thus  Condition S2.1 will be satisfied at the end of Step 2. 
 And as we only ever add good edges to $R_A$ and $R_B$, Condition S2.2 will  hold automatically.
 We  now show that Condition S2.3  will also be satisfied. Recall that 
 $$
 S_A=\{u\in S\cap A\,:\, d_{G_{A,B}^*}(u) \le k-2n^{\frac{2}{3}} \} \quad \text{and}\quad S_B=\{u\in S\cap B\,:\, d_{G_{A,B}^*}(u) \le k-2n^{\frac{2}{3}} \}.
 $$
 Since $\sum_{u\in C\cap S} \df_{G^*_{A,B}}(u)\le k \le \frac{1}{2}\Delta(G)+0.6\eta n<2n$ for any $C\in \{A,B\}$, it follows
 that 
 $$
 |S_A| <n^{1/3} \quad \text{and} \quad |S_B|<n^{1/3}. 
 $$

 In the process of Step 2,   the number of newly colored edges of $H$ that are incident with a vertex $u\in V(G^*)\setminus (S_A\cup S_B\cup N^b(x)\cup \{x\})$ will equal  the number of alternating paths containing $u$ that have been exchanged. 
 The number of such alternating paths of which $u$ is the first vertex will equal the number of colors that missed at $u$ at the end of Step 1.  The number of alternating paths in which $u$ is not the first  vertex will  equal  the degree of $u$ in $R_A\cup R_B$, and so will be less than $r$. Hence the number of colored edges of $H$ that are incident with $u$ will be less than $|\pbar(u)|+r$.  The number of newly colored edges of $H$ that are incident with a vertex $u\in S_A\cup S_B\cup N^b(x)\cup \{x\}$ will equal  the number of colors that missed at $u$ at the end of Step 1. 
  Thus Condition S2.3   will be satisfied.

As the vertex $x$ may have a quite small simple degree, we deal with 
colors in $\pbar(x)$ first. By the upper bound on $|\pbar(x)|$
given in Step 1, we have $|\pbar(x)|<\eta n$ when $G$ is in Conditions~\eqref{Ca} to~\eqref{Cd},
and 
  $|\pbar(x)| \le \frac{1}{2}(1-\ve)n+n^{\frac{2}{3}}$ when $G$ is in Condition~\eqref{Ce}. 
When $G$ is in Condition~\eqref{Ce}, as $\mu(G) \le \frac{2}{\eta}$ and $e_G(x,B)-e_{G_{A,B}^*}(x,B) \ge \frac{n}{2}$,  
we know that $x$ is adjacent in $G$ with at least 
$ \frac{n/2}{2/\eta} \ge \frac{\eta }{4}n$ vertices from $B$
that each joins $x$ with an uncolored edge of $H$.

 For each $i\in \pbar(x)$, if there is a vertex  $w\in N_G(x, B)$ with $xw$  uncolored and $i\in \pbar(w)$, we simply color $xw$ by the color $i$. 
 Thus  we assume that for  every $w\in N_G(x, B)$ with $xw$  uncolored, it holds that $i\not\in \pbar(w)$. We choose $w\in N_G(x, B)$ with $xw$ uncolored and with $d_{R_B}(w)$ smallest. 
 Let $w^*\in N_G(w, B)$ such that $ww^*$ is colored by $i$ (the vertex $w^*$ exists as we only added colored edges between $x$ and $N_G(x, B)$ so far and $i\in \varphi(w)$). 
 Then we color $xw$ by $i$ and uncolor $ww^*$.   (We say $w^*$ is \emph{originated} from $x$ with respect to the color $i$). 
  This procedure 
 guarantees all edges of $R_B$ to be good, as 
 when $G$ is in Conditions~\eqref{Ca} to~\eqref{Cd}, we have 
 $ |\pbar(x)|<\eta n<  \eta^{\frac{1}{2}}n$; and when $G$ is in Condition~\eqref{Ce},
 although $|\pbar(x)| \le \frac{1}{2}(1-\ve)n+n^{\frac{2}{3}}$, 
 we have at least $\frac{1 }{4}\eta n-|\pbar(x)|/n^{\frac{5}{6}}>1$ choices 
 each time to pick a desired vertex $w\in N_G(x, B)$. 
 After this procedure, all the $k$ colors present at $x$. 
 Thus $x$ is not contained in any MCC-pair with respect to the current edge coloring.

 Let $S^*_A=S_A \cup \{x\}$ and  $S_B^*=S_B\cup N^b(x)$. 
 We now show below the existence of alternating paths for the current MCC-pairs.
    For a given color $i\in [1,k]$,  and vertices $a\in A$ and $b\in B$, let $N_B(a)$ be the set
 of vertices in $B$ that are joined with $a$ by an uncolored edge and are incident with a
 good edge colored $i$ such that the good edge is not incident with any vertex of $S^*_B$, and let $N_A(b)$ be the set of vertices in $A$ that are joined with
 $b$ by an uncolored edge and are incident with a good edge colored $i$ such that the good edge is not incident with any vertex of $S^*_A$. 
 In order to estimate the sizes of $N_A(b)$ and $N_B(b)$, we 
 show that $A$ and $B$ contain only a few vertices 
 that either miss the color $i$ 
 or are incident  with a non-good edge colored $i$. 
 By S2.1, there are at most  $4s$ edges in $R_B$, so there are fewer than $\frac{8s}{r}$ vertices of
 degree at least $r$ in $R_B$. 
 Each non-good edge is incident with one or two 
 vertices of $R_B$ through the color $i$, 
 so there are fewer than
 $
 \frac{16s}{r}
 $ 
 vertices in $B$ that are
 incident with a non-good edge colored $i$. 
 Furthermore, there are at most $2|S^*_B| \le 3n^{\frac{1}{3}}$ vertices in $B$
that are  either contained in $S^*_B$ 
 or adjacent to a vertex from $S^*_B$ through an edge with color $i$. 
Finally, there are fewer than $\max\{3\eta n, 7n^{\frac{2}{3}} \}=3\eta n$
 vertices in $B$ that are missed by the color $i$ by S1.2. So the number of vertices in $B$ that are
 not incident with a good edge colored $i$  such that the good edge is not 
 incident with any vertex from $S^*_B$ is less than
 $$
 \frac{16s}{r}+3n^{\frac{1}{3}}+3\eta n<49 \eta^{\frac{1}{2}}n. 
 $$
 Similarly, 
 the number of vertices  in $A$ that are
 not incident with a good edge colored $i$ such that the good edge is not 
 incident with any vertex from $S^*_A$ is less than
 $49 \eta^{\frac{1}{2}}n$. Note that, by the partition $\{A,B\}$ of $V(G)$, we have 
 \begin{eqnarray}
| N_G(v,B)| &\ge& \frac{1}{2}(d^s_G(v)-n^{\frac{2}{3}}-2n^{\frac{1}{2}}- n^{\frac{1}{2}}) \quad \text{for $v\in A\setminus \{x\}$};  \nonumber \\
 &>& \frac{1}{2}(d^s_G(v)-2n^{\frac{2}{3}}) \nonumber \\
 && \label{eqn:v_H-degree} \\
 | N_G(v,A)| &\ge& \frac{1}{2}(d^s_G(v)-n^{\frac{2}{3}}-2n^{\frac{1}{2}}- n^{\frac{1}{2}}) \quad \text{ for $v\in B$}. \nonumber\\ 
 &>& \frac{1}{2}(d^s_G(v)-2n^{\frac{2}{3}}) \nonumber
  \end{eqnarray}
By S2.3 and the upper bound on $|\pbar(u)|$ for any $u\in V(G^*)$ provided in Step 1,
we have the following lower bounds on $|N_B(a)|$ and $|N_A(b)|$: 
\begin{enumerate}[]
	\item When $\{a,b\}\cap (S^*_A\cup S^*_B) = \emptyset$, 
	\begin{equation}\label{eqn: not-in-S-AB}
		|N_A(b)|, |N_B(a)| \ge \frac{1}{2}\left((1+\frac{1}{2}\ve)n-2n^{\frac{2}{3}}\right)-(2\eta n+\eta^{\frac{1}{2}}n)-49\eta^{\frac{1}{2}}n>(\frac{1}{2}+\frac{1}{5}\ve)n. 
	\end{equation} 
\item When $G$ is in Conditions~\eqref{Cb} or~\eqref{Cc} and  $\{a,b\}\cap N^b(x) \ne \emptyset$, 
 we have 
\begin{equation}\label{eqn:condition-b-c}
	|N_A(b)|, |N_B(a)| \ge \frac{1}{2}\left((1+\ve)n-2n^{\frac{2}{3}} \right)-(\frac{1}{2}-\frac{1}{3}\ve)n-49\eta^{\frac{1}{2}}n>\frac{1}{2}\ve n. 
\end{equation}
\item When $G$ is in Condition~\eqref{Cc} and $b=y$, 
 we have 
\begin{equation}\label{eqn:b=y}
	|N_A(b)| \ge \frac{1}{2}\left(2\ve n-2n^{\frac{2}{3}} \right)-\eta n-49\eta^{\frac{1}{2}}n>\frac{3}{4}\ve n. 
\end{equation}
\item 
 When $G$ is in Condition~\eqref{Cd} and $\{a,b\}\cap \{y,z\} \ne \emptyset$,  we have 
\begin{equation}\label{eqn5}
	|N_A(b)|, |N_B(a)| \ge \frac{1}{2}\left((\frac{1}{2}+\frac{3\ve}{2})n-\eta n-2n^{\frac{2}{3}}\right)-(\frac{1}{4}-\frac{1}{5}\ve)n-49\eta^{\frac{1}{2}}n>\frac{3}{4}\ve n. 
\end{equation}
 When $G$ is in Condition~\eqref{Ce} and $\{a,b\}\cap (S^*_A\cup S^*_B) \ne \emptyset$,  we have 
\begin{equation}\label{eqn5}
	|N_A(b)|, |N_B(a)| \ge \frac{1}{2}\left((1+\ve)n-2n^{\frac{2}{3}}\right)-(\frac{1}{2}-\frac{1}{3}\ve)n-49\eta^{\frac{1}{2}}n>\frac{1}{2}\ve n. 
\end{equation}

\end{enumerate}

 Let $M_B(a)$ be the set of vertices in $B$ that are joined with a vertex in $N_B(a)$ by an edge
 of color $i$, and let $M_A(b)$ be the set of vertices in $A$ that are joined with a vertex in $N_A(b)$
 by an edge of color $i$.  Note that $(S^*_A\cup S^*_B)\cap ( M_A(b) \cup M_B(a))=\emptyset$ by the choice of $N_A(b)$ and $N_B(a)$. 
 Note also that $|M_B(a)|=|N_B(a)|$     
 but some vertices
 may be in both. Similarly
 $|M_A(b)|=|N_A(b)|$. 
 
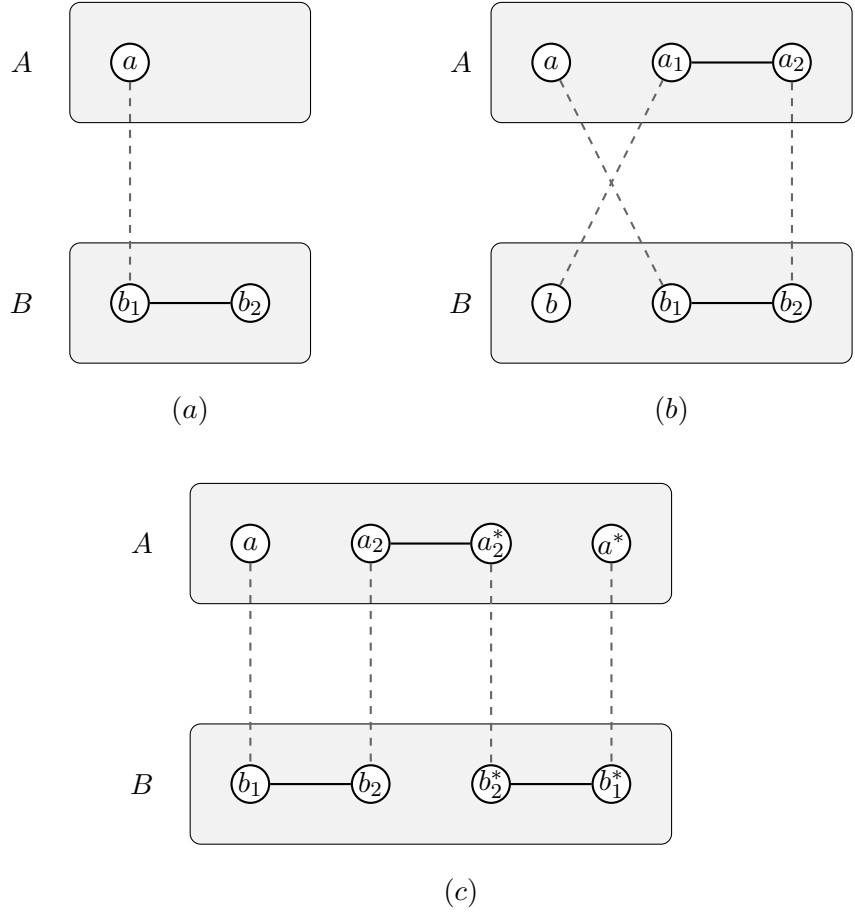
\begin{figure}[!htb]
	\begin{center}
		
		\begin{tikzpicture}[scale=0.8]
		
			\begin{scope}[shift={(-10,0)}]
		\draw[rounded corners, fill=white!90!gray] (8, 0) rectangle (12, 2) {};
		
		\draw[rounded corners, fill=white!90!gray] (8, -4) rectangle (12, -2) {};
		
		{\tikzstyle{every node}=[draw ,circle,fill=white, minimum size=0.5cm,
			inner sep=0pt]
			\draw[black,thick](9,1) node (c)  {$a$};
		%	\draw[black,thick](11,1) node (c1)  {$a_{2}$};
			%\draw[black,thick](12,1) node (c2)  {$a_{12}$};
			%\draw[black,thick](13.5,1) node (c3)  {$a_3$};
			%\draw[black,thick](15,1) node (c4)  {$a_2$};
			\draw[black,thick](9,-3) node (d)  {$b_{1}$};
			\draw[black,thick](11,-3) node (d1)  {$b_{2}$};
			%\draw[black,thick](12,-3) node (d2)  {$b_2$};
			%\draw[black,thick](13.5,-3) node (d3)  {$b_{22}$};
			%\draw[black,thick](15,-3) node (d4)  {$b_{21}$};
			
		}
		\path[draw,thick,black!60!white,dashed]
		(c) edge node[name=la,pos=0.7, above] {\color{blue} } (d)
	%	(c1) edge node[name=la,pos=0.7, above] {\color{blue} } (d1)	
		;
		
		\path[draw,thick,black]
		(d) edge node[name=la,pos=0.7, above] {\color{blue} } (d1)
		;
		\node at (7.2,1) {$A$};
		\node at (7.2,-3) {$B$};
		\node at (10,-4.8) {$(a)$};	
		\end{scope}	
		
		\begin{scope}[shift={(5,0)}]
		\draw[rounded corners, fill=white!90!gray] (0, 0) rectangle (6, 2) {};
		
		\draw[rounded corners, fill=white!90!gray] (0, -4) rectangle (6, -2) {};
		
		{\tikzstyle{every node}=[draw ,circle,fill=white, minimum size=0.5cm,
			inner sep=0pt]
			\draw[black,thick](1,1) node (a)  {$a$};
			\draw[black,thick](3,1) node (a1)  {$a_1$};
			\draw[black,thick](5,1) node (a2)  {$a_2$};
			\draw[black,thick](1,-3) node (b)  {$b$};
			\draw[black,thick](3,-3) node (b1)  {$b_1$};
			\draw[black,thick](5,-3) node (b2)  {$b_2$};
		}

		\path[draw,thick,black!60!white,dashed]
		(a) edge node[name=la,pos=0.7, above] {\color{blue} } (b1)
		(a2) edge node[name=la,pos=0.7, above] {\color{blue} } (b2)
		(b) edge node[name=la,pos=0.6,above] {\color{blue}  } (a1)
		;
		
		\path[draw,thick,black]
		(a1) edge node[name=la,pos=0.7, above] {\color{blue} } (a2)
		(b1) edge node[name=la,pos=0.7, above] {\color{blue} } (b2)
		;
		
		\node at (-0.5,1) {$A$};
		\node at (-0.5,-3) {$B$};
		\node at (3,-4.8) {$(b)$};

		\end{scope}

		\begin{scope}[shift={(-8,-8)}]
		\draw[rounded corners, fill=white!90!gray] (8, 0) rectangle (16, 2) {};
		
		\draw[rounded corners, fill=white!90!gray] (8, -4) rectangle (16, -2) {};
		
		{\tikzstyle{every node}=[draw ,circle,fill=white, minimum size=0.5cm,
			inner sep=0pt]
			\draw[black,thick](9,1) node (c)  {$a$};
			\draw[black,thick](11,1) node (c1)  {$a_{2}$};
			\draw[black,thick](13,1) node (c2)  {$a^*_{2}$};
			\draw[black,thick](15,1) node (c3)  {$a^*$};
			%\draw[black,thick](12,1) node (c2)  {$a_{12}$};
			%\draw[black,thick](13.5,1) node (c3)  {$a_3$};
			%\draw[black,thick](15,1) node (c4)  {$a_2$};
			\draw[black,thick](9,-3) node (d)  {$b_{1}$};
			\draw[black,thick](11,-3) node (d1)  {$b_{2}$};
			\draw[black,thick](13,-3) node (d2)  {$b^*_{2}$};
			\draw[black,thick](15,-3) node (d3)  {$b^*_{1}$};
			%\draw[black,thick](12,-3) node (d2)  {$b_2$};
			%\draw[black,thick](13.5,-3) node (d3)  {$b_{22}$};
			%\draw[black,thick](15,-3) node (d4)  {$b_{21}$};
			
		}
		\path[draw,thick,black!60!white,dashed]
		(c) edge node[name=la,pos=0.7, above] {\color{blue} } (d)
		(c1) edge node[name=la,pos=0.7, above] {\color{blue} } (d1)	
		(c3) edge node[name=la,pos=0.7, above] {\color{blue} } (d3)	
		(c2) edge node[name=la,pos=0.7, above] {\color{blue} } (d2)	
		;
		
		\path[draw,thick,black]
		(d) edge node[name=la,pos=0.7, above] {\color{blue} } (d1)
		(c2) edge node[name=la,pos=0.7, above] {\color{blue} } (c1)
		(d3) edge node[name=la,pos=0.7, above] {\color{blue} } (d2)
		;
		\node at (7.2,1) {$A$};
		\node at (7.2,-3) {$B$};
		\node at (12.5,-4.8) {$(c)$};	
		\end{scope}	
		
		\end{tikzpicture}
		-	  	\end{center}
	\caption{The alternating path $P$. Dashed lines indicate uncoloured edges, and solid
		lines indicate edges with color $i$.}
	\label{f1}
\end{figure}

 For a MCC-pair $(a,b)$, in order to have a unified discussion as in the case that $\{a,b\}\cap (S^*_A\cup S^*_B) = \emptyset$,  if necessary, by exchanging an alternating path of length 2 from $a$ to another vertex $a^*$,
 and exchanging an alternating path from $b$ to another vertex $b^*$, 
 we will replace the pair $(a,b)$ by $(a^*,b^*)$ such that $\{a^*,b^*\}\cap (S^*_A\cup S^*_B)=\emptyset$.  
 Precisely, we will implement the following 
 operations to vertices in $S_A\cup S^*_B$. 
 For any vertex $a\in S_A$, and for each color $i\in \pbar(a)$, we take an edge $b_1b_2$ with $b_1\in N_B(a)$ and $b_2\in M_B(a)$ such that $b_1b_2$ is colored by $i$, where the edge $b_1b_2$ exists by~\eqref{eqn: not-in-S-AB}-\eqref{eqn5} and the fact that $|M_B(a)|=|N_B(a)|$.   Then we exchange the path $ab_1b_2$ by coloring $ab_1$
 with $i$ and uncoloring the edge $b_1b_2$ (See Figure~\ref{f1}(a)).
  After this, the edge $ab_1$  of $H$ is now colored by $i$, and the uncolored edge $b_1b_2$ is added to $R_B$.
 We then update the original MCC-pair that contains $a$ with respect to the color $i$ by replacing the vertex $a$ with $b_2$. 
 We do this at the vertex $a$ for every color $i\in \pbar(a)$ and then repeat the same process for every vertex in $S_A$. 
 Similarly, 
 for any vertex $b\in S^*_B$, and for each color $i\in \pbar(b)$, we take an edge $a_1a_2$ with $a_1\in N_A(b)$ and $a_2\in M_A(b)$ such that $a_1a_2$ is colored by $i$, where the edge $a_1a_2$ exists by~\eqref{eqn: not-in-S-AB}-\eqref{eqn5} and the fact that $|M_B(b)|=|N_B(b)|$.   Then we exchange the path $ba_1a_2$ by coloring $ba_1$
 with $i$ and uncoloring the edge $a_1a_2$. The same, we update the original MCC-pair that contains $b$ with respect to the color $i$ by replacing the vertex $b$ with $a_2$.

After the procedure above,  we have now three types MCC-pair $(u,v)$: $u,v\in A$, $u,v\in B$, and $A$
contains exactly one of $u$ and $v$ and $B$ contains the other. However, in either case, $\{u,v\}\cap (S^*_A\cup S^*_B)=\emptyset$. We will exchange alternating path for each of such pairs. 

We deal with each of the colors from $[1,k]$ in turn. 
Let $i\in [1,k]$ be a color.  We consider first an MCC-pair $(a,a^*)$
with respect to  $i$ such that $a,a^* \in A$. 
By~\eqref{eqn: not-in-S-AB}, we have  $|M_B(a^*)|>(\frac{1}{2}+\frac{1}{5}\ve)n$.
We take an edge $b_1^*b_2^*$ colored by $i$ with $b_1^* \in N_B(a^*)$ and  $b_2^*\in M_B(a^*)$.  
Then again, by~\eqref{eqn: not-in-S-AB}, we have 
$
 |M_B(a)|,  |M_A(b_2^*)|>(\frac{1}{2}+\frac{1}{5}\ve)n.  
$
Therefore, as each vertex  $c\in M_A(b_2^*)$ 
satisfies $|N_B(c)|>(\frac{1}{2}+\frac{1}{5}\ve)n$,  we have $|N_B(c)\cap M_B(a)| \ge \frac{2}{5}\ve n$. We take $a_2a_2^*$ colored by $i$ with $ a_2^*\in N_A(b_2^*)$ and  $a_2\in   M_A(b_2^*)$. 
Then we let $b_2\in N_B(a_2)\cap M_B(a)$, and let $b_1$
be the vertex in $N_B(a)$ such that $b_1b_2$ is colored by $i$. 
Now we get the alternating path $P=ab_1b_2 a_2 a_2^* b_2^* b_1^* a^*$ (See Figure~\ref{f1}(c)).  
We exchange $P$ by coloring $ab_1, b_2a_2, a_2^*b_2^*$ and $b_1^*a^*$
with color $i$ and uncoloring the edges $b_1b_2, b_1^*b_2^*$ and $a_2a_2^*$. 
After the exchange,  the color $i$ appears on edges incident with $a$ and $a^*$,
the edges $b_1b_2$ and $b_1^*b_2^*$ are added to $R_B$
and the edge $a_2a_2^*$ is added to $R_A$. We added at most 
one edges to each of $R_A$ and $R_B$ when we updated the original MCC-pair corresponding to $(a,a^*)$, where $a$ and $a^*$ could be originated from some vertices $b_1', b_2'\in S_B^*$. 
Also, exactly one of $b_1'$ and $b_2'$ could be originated from  $x$
with respect to the color $i$.    
Thus  we added at most three edges to $R_A$
and at most four edges to $R_B$ when we modify $\varphi$ to have the color $i$ present at both of the  vertices in the original MCC-pair corresponding to $(a,a^*)$.  The maximum length of an alternating path combined from the three procedures (dealing missing colors at $x$, dealing missing colors at vertices from $S_A\cup S_B^*$, and dealing missing colors at vertices from $V(G^*)\setminus (S_A^*\cup S_B^*)$) together is at most $2+2+2+7=13$. 
By symmetry, we can  deal with an MCC-pair $(b,b^*)$
with respect to  $i$ such that $b,b^* \in B$ similarly as above.

Then we consider an MCC-pair $(a,b)$
with respect to  $i$ such that $a\in A$ and $b\in B$. 
By~\eqref{eqn: not-in-S-AB}, we have 
$
|M_B(a)|,  |M_A(b)|>(\frac{1}{2}+\frac{1}{5}\ve)n.  
$
We choose $a_1a_2$ with color $i$ such that $a_1\in N_A(b)$
and $a_2\in M_A(b)$. Now as $|M_B(a)|, |N_B(a_2)| >(\frac{1}{2}+\frac{1}{5}\ve)n$ by~\eqref{eqn: not-in-S-AB}, we know that $N_B(a_2)\cap M_B(a) \ne \emptyset$. We choose $b_2\in N_B(a_2)\cap M_B(a)$
and let $b_1\in N_B(a)$ such that $b_1b_2$ is colored by $i$. 
Then $P=ab_1b_2a_2a_1 b$ is an alternating path from $a$ to $b$ (See Figure~\ref{f1}(b)). We exchange $P$ by coloring $ab_1, b_2a_2$ and $a_1b$
with color $i$ and uncoloring the edges $a_1a_2$ and $b_1b_2$. 
After the exchange,  the color $i$ appears on edges incident with $a$ and $b$,
the edge $a_1a_2$ is added to $R_A$
and the edge $b_1b_2$ is added to $R_B$. We added at most 
one edge to each of $R_A$ and $R_B$ when we updated the original MCC-pair corresponding to $(a,b)$. Thus  we added at most one edge to $R_A$ and at most two edges to $R_B$ when we modify $\varphi$ to have the color $i$ present at both of the vertices in the original MCC-pair corresponding to $(a,b)$. 
 By finding
such paths  for all MCC-pairs with respect to the color $i$,  we can increase the number of
edges colored $i$ until the color class is a 1-factor of $G^*$. By doing this for all colors,
we can make each of the $k$ color classes a 1-factor of $G^*$.

\begin{center}
	Step 3: Coloring $R_A$ and $R_B$ and extending the new color classes
\end{center}

Each of the color classes for the colors from $[1,k]$ is now a 1-factor of $G^*$. We 
now consider the multigraphs $R_A$ and $R_B$ that consist of the uncolored edges of $G^*_A$ and $G^*_B$. 
By Condition S2.1, $R_A$ and $R_B$ each has fewer than $4s$ edges, and $\Delta(R_A), \Delta(R_B)< r$.  Note that $R_A$
and $R_B$ might contain parallel edges 
with endvertices in $S$ or between $y$ and $z$ when $G$ is in Condition~\eqref{Cd}. 
 By Theorem~\ref{chromatic-index} and  Theorem~\ref{lem:equa-edge-coloring},  $R_A$ and $R_B$ each has an
equalized edge-coloring with exactly $\ell:= 2\lfloor r \rfloor$ colors $k+1, \ldots, k+\ell$. 

If $G$ is in Conditions~\eqref{Ca} to~\eqref{Cd},   then we have $e(R_A)=e(R_B)$. 
Under these  conditions, by renaming some color classes of $R_A$ if necessary,  we can assume that in the edge colorings of $R_A$ and $R_B$, each color
appears on the same number of edges in $R_A$ as it does in $R_B$. 
When $G$ is in Condition~\eqref{Ce}, by our assumption that $G^*_B$
has more edges than $G^*_A$ does, we have  $e(R_A) \le e(R_B)$. 
In this case, we can assume that in the edge colorings of $R_A$ and $R_B$, the number of edges with a color $i\in [k+1,k+\ell]$ in $R_B$
is at least the number of edges with a color $i\in [k+1,k+\ell]$ in $R_A$.

There are fewer than $4s$ edges in each of $R_A$ and $R_B$, and $\ell > r$, so each of
the color $i\in [k+1,k+\ell]$ appears on fewer than
$\frac{4s}{r}$
edges in each of $R_A$ and $R_B$. We will now color some of the edges of $H$ with the $\ell$
colors from $[k+1,k+\ell]$ so that each of these color classes present at 
vertices from $V(G^*)\setminus U$. 
We  perform the following procedure for each of the $\ell$ colors in turn.

 Given a
color $i$ with $i\in [k+1,k+\ell]$, we let $A_i$ and $B_i$ be the sets of vertices in $A$ and
$B$ respectively that are incident with edges colored $i$. 
Note that $|A_i| \le  |B_i| <2\times \frac{4s}{r} =\frac{8s}{r}$, 
 as $R_A$ and $R_B$ each contains fewer than $ \frac{4s}{r}$ edges colored $i$.   Note that if $G$ is in Conditions~\eqref{Ca} to~\eqref{Cd}, we have $|A_i|=|B_i|$;
 and we might have $|B_i| > |A_i|$ when $G$ is in Condition~\eqref{Ce}. 
 When $G$ is in Condition~\eqref{Ce} and $|B_i| > |A_i|$, we let 
 $$
 A_i^* \subseteq (U\cap A) \setminus A_i 
 $$
 such that $|A_i^*|+|A_i|=|B_i|$, and just let $A_i^*=\emptyset$ otherwise.  
 Note that such $A_i^*$ exists as $|U\cap A| \ge \frac{1}{2}\eta n-1$
 and $|A_i|, |B_i| < \frac{8s}{r}=56n^\frac{5}{6}$. 
Let $H_i$ be the subgraph of $H$ obtained by
deleting the vertex sets $A_i \cup A_i^*$ and $B_i$ and removing all colored edges. We will show next that $H_i$ has a perfect matching and we will color 
the edges in the matching by the color $i$. 

For each $u\in V(G^*)$,  the number of colored edges of $H$ incident with 
$u$ is specified by S2.3 and at
most $\ell$ edges incident with $u$ have been colored in Step 3. Also 
each vertex in $G^*$ has fewer than $\frac{8s}{r} <24 \eta^{\frac{1}{2}} n$ edges that join it with a vertex in $A_i$ or $B_i$. So by~\eqref{eqn:v_H-degree}, each vertex  $u\in V(H_i)\setminus (S^*_A\cup S^*_B)$ satisfies 
$$
d_{H_i}^s(u)>  \frac{1}{2}\left((1+\frac{1}{2}\ve)n-2n^{\frac{2}{3}}\right)-(2\eta n+\eta^{\frac{1}{2}}n)-\ell-24 \eta^{\frac{1}{2}}n>(\frac{1}{2}+\frac{\ve}{5})n. 
$$

Note that $V(R_A\cup R_B)\cap (S_A^*\cup S_B^*)=\emptyset$. 
When $G$ is in Conditions~\eqref{Cb} or~\eqref{Cc}, each of the  vertex  $u\in N^b(x)$, 
where $u\ne y$ when $G$ is in Condition~\eqref{Cc},  satisfies 
$$
d_{H_i}^s(u)>\frac{1}{2}\left((1+\ve)n-2n^{\frac{2}{3}}\right)- (\frac{1}{2}-\frac{1}{3}\ve )n-\ell-24 \eta^{\frac{1}{2}}n >\frac{1}{2} \ve n.
$$

When $G$ is in Condition~\eqref{Cc}, the vertex 
$y$ satisfies 
$$
d_{H_i}^s(y) >\frac{1}{2}\left(2\ve n-2n^{\frac{2}{3}} \right)-\eta n-\ell-24 \eta^{\frac{1}{2}}n>\frac{\ve}{4} \ve n. 
$$

When $G$ is in Condition~\eqref{Cd}, each  $u\in \{y,z\}$ satisfies  
$$ d_{H_i}^s(u)>\frac{1}{2}\left((\frac{1}{2}+\frac{3\ve}{2})n-\eta n-2n^{\frac{2}{3}}\right)-(\frac{1}{4}-\frac{1}{5}\ve)n-\ell-24 \eta^{\frac{1}{2}}n>\frac{3}{4}\ve n.
$$

When $G$ is in Condition~\eqref{Ce}, each $u\in S_A\cup S_B$
satisfies 
$$d_{H_i}^s(u)>\frac{1}{2}\left((1+\ve)n-2n^{\frac{2}{3}}\right)-(\frac{1}{2}-\frac{1}{3}\ve)n-\ell-24\eta^{\frac{1}{2}}n>\frac{1}{2}\ve n.  
$$

For the vertex $x$, we have  
$$
d_{H_i}(x) \ge \delta(G)-k-\ell-24\eta^{\frac{1}{2}}n \ge (1+\ve)n- \left(\frac{1}{2}\Delta(G)+0.6\eta n\right)-\ell -24\eta^{\frac{1}{2}}n>\frac{1}{2}\ve n. 
$$

Thus $\delta(H_i) \ge 1$ and $\delta(H_i-x) \ge \frac{1}{2}\ve n-1$ in any case and 
$H_i$ has at most $|S^*_A\cup S^*_B| \le 3n^{1/3} <\frac{1}{2} \ve n-1$
vertices of simple degree less than $\frac{1}{2}n+1$. 
So $H_i$
has a 1-factor $F$ by Lemma~\ref{lem:matching-in-bipartite}. If we color the edges of $F$ with the color
$i$, then every vertex in $V(G^*)\setminus A_i^*$ is incident with an edge of color $i$. 
We repeat this procedure for each of the colors from $[k+1,k+\ell]$. After this has been
done, each of these $\ell$ colors presents at all vertices from $V(G^*)\setminus U$. So at the conclusion of Step 3,
all of the edges in $G^*_A$ and $G^*_B$ are colored, some of the edges of $H$ are colored, 
each of the $k$ color classes  for colors from $[1,k]$ is a 1-factor of $G^*$,
and  each of the $\ell$ colors  from $[k+1,k+\ell]$ 
presents at all vertices from $V(G^*)\setminus U$.

\begin{center}
	Step 4: Coloring the graph $R$ 
\end{center}
Let $R$ be the submultigraph of $G^*$ consisting of the remaining uncolored edges. These
edges all belong to $H$, so $R$ is a submultigraph of $H$ and hence is bipartite. 
We claim that $\Delta(R) = \Delta(G^*)-k-\ell$. 
This is obvious when $G$ is in Conditions~\eqref{Ca} to~\eqref{Cd}. 
Thus we assume that $G$ is in Condition~\eqref{Ce}. 
Note that every vertex from $V(G^*)\setminus U$ 
presents every color from $[1,k+\ell]$
and so those vertices have degree  at most  $\Delta(G^*)-k-\ell$ 
in $R$. 
For some vertices from $U$, they present all
the colors from $[1,k]$. Thus by~\eqref{eqn:degree-U-vertex-in-G*}, those vertices 
have degree at most  
$$
\Delta(G^*)-\frac{1}{3} \eta n -k <\Delta(G^*)-k-\ell = \Delta(G^*)-k-2\lfloor n^{\frac{5}{6}} \rfloor
$$
in $R$. 
By Theorem~\ref{konig} we can color the
edges of $R$ with $\Delta(R)$ colors from $[k+\ell+1, \Delta(G^*)]$. 
Thus $\chi'(G^*) \le k+\ell +(\Delta(G^*)-k-\ell)=\Delta(G^*)$
and so $\chi'(G^*)=\Delta(G^*)$, as desired. 

Lastly, we check that there is a polynomial time algorithm to obtain an
edge $\Delta(G)$-coloring of $G$.
By Lemma~\ref{lem:partition}, we can obtain a desired partition  
$\{A,B\}$ of $V(G)$ in polynomial time. Also,  it is polynomial time 
to edge color $G_{A,B}$
by an algorithm described in~\cite{MR1156837}. Modifying $G_A$
and $G_B$ into $G_A^*$ and $G_B^*$ and the corresponding edge colorings
into the desired edge-colorings can be done in polynomial time too.  
In Step 2, the construction of the alternating paths  and swaps of the colors on the paths can be done in $O(n^3)$-time, as the total number of colors missing at vertices is $O(n^2)$ and  it takes  $O(n)$-time  to find an alternating path for a MCC-pair.  
In Step 3, there is polynomial time algorithm (see e.g. \cite{MR875324}) to edge color $R_A$ and $R_B$ using at most $\ell$ colors. Then by 
doing Kempe changes as  mentioned in the comments immediately after Theorem~\ref{lem:equa-edge-coloring}, these edge colorings can be modified
into  equalized edge-colorings in polynomial time. 
The last step is to edge color the bipartite multigraph $R$ using $\Delta(R)$ colors, which can be done in polynomial time in $n$, for example, using an algorithm from~\cite{MR664720}. 
Thus, there is a polynomial time 
algorithm that gives an edge coloring of  $G$ using $\Delta(G)$ colors. 
\qed

%\section*{Acknowledgment}

%\bibliographystyle{plain}
%\bibliography{SSL-BIB_08-19}

\begin{thebibliography}{10}
		\bibitem{MR0411988}
	J.~A. Bondy and U.~S.~R. Murty.
	\newblock Graph theory with applications. 
	\newblock{\em American Elsevier Publishing Co., Inc., New York}, 
	(Elsevier North-Holland), 1976. 
	
		\bibitem{MR1001390}
	A.~G. Chetwynd and A.~J.~W. Hilton.
	\newblock {$1$}-factorizing regular graphs of high degree---an improved
	bound. 
	\newblock {\em Discrete Math.}, 75(1-3): 103--112, 1989. 
	
	
	
	\bibitem{MR848854}
	A.~G. Chetwynd and A.~J.~W. Hilton.
	\newblock Star multigraphs with three vertices of maximum degree.
	\newblock {\em Math. Proc. Cambridge Philos. Soc.}, 100(2):303--317, 1986.
	
	
	
	

	\bibitem{MR975994}
	A.~G. Chetwynd and A.~J.~W. Hilton.
	\newblock The edge-chromatic class of graphs with maximum degree at least
	{$|V|-3$}.
	\newblock In {\em Graph theory in memory of {G}. {A}. {D}irac ({S}andbjerg,
		1985)}, volume~41 of {\em Ann. Discrete Math.}, pages 91--110. North-Holland,
	Amsterdam, 1989.
	

	
	\bibitem{MR664720}
	R. Cole and J. Hopcroft.
	\newblock On edge coloring bipartite graphs.
	\newblock {\em SIAM J. Comput.}, 11(3):540--546, 1982.
	
	\bibitem{MR3545109}
	B. Csaba, D. K\"{u}hn, A. Lo, D. Osthus, and A. Treglown.
	\newblock Proof of the 1-factorization and {H}amilton decomposition
	conjectures.
	\newblock {\em Mem. Amer. Math. Soc.}, 244(1154):v+164, 2016.
	
	\bibitem{MR47308}
	G.~A. Dirac.
	\newblock Some theorems on abstract graphs.
	\newblock {\em Proc. London Math. Soc. (3)}, 2:69--81, 1952.
	
	\bibitem{fw}
	S.~Fiorini and R.~J. Wilson.
	\newblock Edge-colourings of graphs, {R}esearch notes in {M}aths.
	\newblock {\em Pitman, London}, 1977.
	
%	\bibitem{GKO}
%	S.  Glock, D. K\"{u}hn, and D.  Osthus.
%	\newblock Optimal path and cycle decompositions of dense quasirandom graphs.
%	\newblock {\em J. Combin. Theory Ser. B}, 118:88--108, 2016.
	
	\bibitem{MR2028248}
	S. Gr\"{u}newald and E. Steffen.
	\newblock Independent sets and 2-factors in edge-chromatic-critical graphs.
	\newblock {\em J. Graph Theory}, 45(2):113--118, 2004.
	
	\bibitem{Gupta-67}
	R.~G. Gupta.
	\newblock {\em Studies in the Theory of Graphs}.
	\newblock PhD thesis, Tata Institute of Fundamental Research, Bombay, 1967.
	
	\bibitem{MR148049}
	S.~L. Hakimi.
	\newblock On realizability of a set of integers as degrees of the vertices of a
	linear graph. {I}.
	\newblock {\em J. Soc. Indust. Appl. Math.}, 10:496--506, 1962.
	
	\bibitem{Holyer}
	I. Holyer.
	\newblock The {NP}-completeness of edge-coloring.
	\newblock {\em SIAM J. Comput.}, 10(4):718--720, 1981.
	
	
		\bibitem{matching}
	J. E. {Hopcroft}, and R. M. {Karp}. 
	\newblock An $n^{5/2}$ algorithm for maximum matchings in bipartite graphs.
	\newblock {\em SIAM J. Comput.}, 2(4):225--231, 1973.
	
	\bibitem{MR875324}
	H.~J. Karloff and D.~B. Shmoys.
	\newblock Efficient parallel algorithms for edge coloring problems.
	\newblock {\em J. Algorithms}, 8(1):39--52, 1987.
	
	\bibitem{MR1511872}
	D.  K\"{o}nig.
	\newblock \"{U}ber {G}raphen und ihre {A}nwendung auf {D}eterminantentheorie
	und {M}engenlehre.
	\newblock {\em Math. Ann.}, 77(4):453--465, 1916.
	
	\bibitem{MR300623}
	C.~J.~H. McDiarmid.
	\newblock The solution of a timetabling problem.
	\newblock {\em J. Inst. Math. Appl.}, 9:23--34, 1972.
	
		
	\bibitem{MR1156837}
	J.~Misra and D.  Gries.
	\newblock A constructive proof of {V}izing's theorem.
	\newblock {\em Inform. Process. Lett.}, 41(3):131--133, 1992.
	

	
	\bibitem{MR1439301}
	L.~Perkovi\'c  and B.~ Reed. 
	\newblock Edge coloring regular graphs of high degree. 
	\newblock {\em Discrete Math.}, 165/166: 567--578, 1997. 
	
	\bibitem{MR2082738}
	M.  Plantholt.
	\newblock Overfull conjecture for graphs with high minimum degree.
	\newblock {\em J. Graph Theory}, 47(2):73--80, 2004.
	
	
	\bibitem{Plantholt2}
	M.  Plantholt.
	\newblock The chromatic index of graphs with large even order $n$ and minimum degree at least $2n/3$.
	\newblock {\em Discrete Math.}, 345(7): 112880, 2022.
	
\bibitem{2104.06253}
M.  Plantholt and S. Shan.
\newblock Edge coloring graphs with large minimum degree.
\newblock {\em \arxiv{2105.05286}}, 2021.


	
	
	
	
	\bibitem{seymour79}
	P.~D. Seymour.
	\newblock On multicolourings of cubic graphs, and conjectures of {F}ulkerson
	and {T}utte.
	\newblock {\em Proc. London Math. Soc. (3)}, 38(3):423--460, 1979.
	
	\bibitem{2104.06253}
	S. Shan.
	\newblock Chromatic index of dense quasirandom graphs.
	\newblock {\em \arxiv{2104.06253}}, 2021.
	
	\bibitem{StiebSTF-Book}
	M.~Stiebitz, D.~Scheide, B.~Toft, and L.~M. Favrholdt.
	\newblock {Graph Edge Coloring: Vizing's Theorem and Goldberg's Conjecture.}
	\newblock {\em Wiley Series in Discrete Mathematics and Optimization}. John Wiley \&
	Sons, Inc., Hoboken, NJ, 2012.
	\newblock 
	
	
	\bibitem{Vizing-2-classes}
	V.~G. Vizing.
	\newblock Critical graphs with given chromatic class.
	\newblock {\em Diskret. Analiz No.}, 5:9--17, 1965.
	
\end{thebibliography}

\end{document}